\documentclass[11pt,a4paper]{article}

\usepackage{graphicx} 
\graphicspath{{Figures/}{figures/}}
\usepackage{subcaption}
\usepackage[a4paper,top=3cm,bottom=3cm,left=3cm, right=3cm]{geometry}
\usepackage{xcolor}
\usepackage{amsmath,bm}
\usepackage{amssymb,amsthm}
\usepackage{orcidlink}
\usepackage{natbib}
\usepackage[noblocks]{authblk}
\usepackage[markup]{changes}
\definechangesauthor[name={Nestor}, color=orange]{NR}
\definechangesauthor[name={Ignacio}, color=teal]{IR}
\newcommand{\nablabar}{\overline{\nabla}}

\theoremstyle{remark}

\newtheorem{remark}{Remark}

\title{A stabilized, nodal strain finite element method \\ 
for small and large deformations:\\
  formulation, analysis, and impact on remeshing strategies}
\author[1]{Nestor Rossi}
\author[2,1]{Ignacio Romero\orcidlink{0000-0003-0364-6969}}

\affil[1]{IMDEA Materials Institute, Eric Kandel 2, 28096 Getafe, Madrid, Spain}
\affil[2]{Dept. of Mechanical Engineering, Universidad Politécnica de Madrid,
  Jos\'{e} Guti\'{e}rrez Abascal, 2, 28006 Madrid, Spain}

\begin{document}

\maketitle

\begin{abstract}
  We present in this article the discretization of small and finite strain mechanics using
  variationally consistent, stabilized, nodal strain finite element formulations. We prove that
  these methods derive from a mixed variational principle that uses different meshes for the primal
  and dual variables. We demonstrate that, for both mechanical problems, nodal strain finite
  elements can be written as pure primal formulations with assumed strain operators. The
  stabilization terms are shown to be necessary and a new, general class of stabilizing functions is
  proposed for the finite strain regime. The ensuing formulations can be interpreted as pure
  particle methods because all kinematic and material information is stored at the nodes. This point
  of view explains the favorable properties of nodal strain methods when employed in problems that
  require frequent remeshing and involve inelastic materials. In these cases, the diffusion of
  internal variables is kept to a minimum. The numerical examples presented validate the claims.
\end{abstract}

\section{Introduction}


The simulation of solids undergoing large deformations is central to computational mechanics. The
complexity of these simulations is often aggravated when solids involve materials modelled with inelastic
constitutive laws. These laws, in addition to being more elaborate than elastic ones, often involve
tracking internal variables --- such as plastic strain or damage --- which need to be stored at
material points (see, e.g., \cite{desouzaneto2011ve,lemaitre2005uh,doghri2000tn}). The accuracy of the overall simulations depends, therefore, not
only on the discretization method employed, but also on the ability to represent faithfully the
history of the material points.


The finite element method is the most widely used numerical approach to solve problems in nonlinear
solid mechanics~\cite{simo2000vm}. It is well known that the accuracy of these approximations strongly
depends on the quality of the computational mesh \cite{babuvska1976angle, lee1993effects}. However, in problems involving
large deformations, Lagrangian finite elements may become severely distorted. Such an unavoidable
process leads to a loss of accuracy in the approximation and, sometimes, even to the appearance of
non-physical failure modes. To alleviate these situations, remeshing techniques are commonly
employed to improve the mesh quality during the simulation. These procedures typically involve
operations such as mesh smoothing, node relocation, local refinement or coarsening, and the reconstruction of element connectivity~\cite{freitag1997tetrahedral}.

Remeshing, however, has undesired side effects. Whenever a mesh is modified, the finite element solution
must be transferred from the discarded mesh to the new one. This transfer involves not only the
primary unknown fields, but also the material history stored in the internal variables, if any. Broadly
speaking, the transfer procedure requires, first, projecting nodal quantities (displacements,
velocities, etc.) from the old mesh to the new mesh. Second, for the internal variables stored at
integration points, an additional operation is required: the information is usually projected first
from the old integration points to the old nodes, and then, once the new mesh is constructed, from
the new nodes to the new integration points. This operation, or its variations, is commonly referred to as history \textit{remapping}.

In standard finite element formulations, even when the nodal positions remain unchanged and only the
element connectivity is modified, the location of the integration points generally changes and thus
a remapping of internal variables is still necessary. This operation introduces numerical diffusion
in the material data and, therefore, loss of information, especially when the internal variables exhibit localized spatial distributions~\cite{mediavilla2006robust}. Such errors directly affect the accuracy of history-dependent simulations, and become especially relevant in problems that require frequent remeshing.


Motivated by its impact on finite strain simulations, a possible strategy to minimize the damaging effects of numerical diffusion during history remapping is to reuse nodes not only to carry the primal variables of the
problem, but also the material history. Naturally, the standard (primal) finite
element method cannot be used for this purpose, since the isoparametric finite element functions do not
have well-defined gradients at the nodes and fluxes (heat flux, strain, etc.) cannot be
calculated at those points.

When the material data is associated with the nodes and remeshing is
performed by exclusively modifying the element connectivity, 
there is no need to transfer internal variables between meshes, and the diffusion of internal
variables is completely bypassed. This possibility has been noted in different contexts
\cite{puso2008meshfree, yuan2019development}, and the idea of combining remeshing with nodal strain
methods has been proposed for the particle finite element method (PFEM) \cite{zhang2018smoothed}. This method --- originally proposed for solving flows with free surfaces \cite{idelsohn2004particle} --- builds on treating nodes as material
particles during remeshing operations~\cite{cremonesi2020state}. However, since it is based on the standard
Gaussian quadrature, remapping procedures are still needed, leading to the diffusion of the internal variables, even if element connectivity remains unchanged. 
There have been some efforts to reduce numerical diffusion during the remapping in the context of PFEM~\cite{rodriguez2016particle, rodriguez2017continuous}.  
We note, in general, though, that the effects of keeping the material information at the nodes on the diffusion of internal variables for history-dependent solid mechanics problems have not yet been systematically quantified.


Several numerical methods based on nodal strain or nodal integration concepts have been proposed in the literature. These approaches often rely on averaging strain fields over nodal patches, which naturally leads to nodal representations of stresses and state variables. Nodal integration techniques have received considerable attention, partly because of their ability to alleviate locking in simplicial meshes \cite{bonet1998simple,dohrmann2000node, puso2006stabilized, krysl2008locking, broccardo2009assumed, liu2009node}. Although this is an important advantage, a detailed assessment of their locking behaviour lies outside the scope of the present work; we refer the reader to the extensive studies cited above. Despite their success, many of these formulations share limitations and unresolved theoretical issues.

A first difficulty is that often the strain-averaging procedure is introduced in a
largely heuristic manner. This issue appears across a variety of closely related approaches reported
under different names in the literature such as nodal averaging methods~\cite{bonet1998simple,dohrmann2000node}, the NICE formulation~\cite{krysl2008locking,
  broccardo2009assumed}, and methods based on assumed mollifiers or smoothing functions, including
averaged nodal strain methods \cite{bonet2001averaged}, NSFEM \cite{liu2009node, nguyen2010node},
LC-PIM \cite{liu2005linearly}, and related meshless techniques \cite{chen2001stabilized}.

A second difficulty concerns the stability of these methods and the stabilization strategies
proposed to address it. The presence of non-physical low-energy deformation modes has been
recognized in dynamic problems \cite{dohrmann2000node, liu2008generalized} and eigenvalue analyses
\cite{puso2006stabilized, broccardo2009assumed}. Moreover, it has been observed that some nodal
integration formulations lose stability as the mesh is refined, due to the lack of coercivity of the
discrete bilinear form \cite{puso2006stabilized, artioli2014assumed}. Nevertheless, stabilization
measures are often omitted in static analyses due to the apparent reliability of nodal approaches,
particularly on coarse meshes \cite{sivapuram2019energy}.

Along with the lack of stability, several stabilization mechanisms have been proposed to overcome
these difficulties, although in most cases without a supporting mathematical analysis. These include
penalizing the difference between nodal and displacement-compatible strain fields through the
introduction of a second material response and an empirically chosen stabilization coefficient
\cite{puso2006stabilized, gee2009uniform, broccardo2009assumed, sivapuram2019energy}; penalizing the
residual of the strong form of the balance of linear momentum \cite{zhang2010temporal,
  feng2014temporal}; and modifying the bilinear form through changes in the smoothing domain,
together with assumptions on the continuity and differentiability of the smoothed strain field
\cite{tong2018high, chen2019gradient}. As a result, the theoretical foundations of these methods, as
well as the consistency of their stabilization terms, are not always clearly
established.
The work of Lamichhane~\cite{lamichhane2009inf, lamichhane2009hu, lamichhane2014mixed}, however, significantly clarified the mathematical foundations of methods based on
nodal averaging. Of particular
relevance to the present work, the averaged nodal strain method was derived consistently
from a Hu-Washizu variational principle~\cite{lamichhane2009hu}, rather than being introduced
as an \emph{ad hoc} numerical assumption. The results of this work were restricted to small strain linear elasticity, leaving open the question of how this framework can be extended to more general nonlinear problems.



The present work builds upon Lamichhane's framework, extending it to large-strain kinematics, and
to constitutive models with internal variables. In the proposed formulation, the strain field and
the history variables of the constitutive model are represented as piecewise constant quantities
defined over nodal patches, which form a dual partition of the domain. As part of this extension,
the projection or averaging operation is adapted to account for changes in the reference
configuration. In addition, a general structure to introduce stabilization as part of the variational
setting is presented which allows greater flexibility when defining the stabilizing function.

The resulting nodal structure described in this work admits a twofold interpretation. On the one hand, the method can be
understood as a mixed finite element method. On the other
hand, it can be interpreted as a particle-like method in which the material information is stored
at, or naturally associated with, the nodes, which behave as material points. This interpretation
makes the formulation particularly well suited for remeshing strategies in which the element
connectivity is modified while the nodal material information is preserved, as discussed before. 
The article investigates this feature and shows that the proposed formulation can substantially reduce the loss of information due to remeshing diffusion, while maintaining a stabilization mechanism consistent with the underlying variational framework.

As a result of the foregoing ideas, the current work is the first to propose a finite strain,
nodally integrated, consistently stabilized finite element formulation for solid mechanics supported by a variational framework. Based on
a rigorous stabilization result in the small strain regime, a stabilization is proposed for the
large strain case that preserves objectivity. The numerical examples shown prove that the proposed
method is robust for finite strain formulations, and reduces the artificial diffusion that results from
remeshing. 

The paper is organized as follows. Section \ref{sect-mixedform-small} reviews the variational
formulation of the nodal averaging approach in the context of small strain kinematics and discusses
some misconceptions found in the literature. In particular, Section \ref{sect-review-stabilization}
illustrates that stabilization is necessary even for static problems on coarse meshes. Section
\ref{sect-mixedform-large} extends the mathematical framework to large-strain kinematics, and
generalizes the stabilization function from the small strain case and the projection operator to account for
changes of reference configuration. Section \ref{sect-remesh-diff} presents the remeshing strategy
and the adopted measure of diffusion, while Section \ref{sect-examples} reports the numerical
results. Finally, Section \ref{sect-conclusion} summarizes the main conclusions of this work. Appendix \ref{sect-appendix} contains information related to the performance of the remeshing algorithm.

\section{Mixed formulation and nodal strain approximation}
\label{sect-mixedform-small}

This section introduces the variational foundation of the nodal strain finite element method for
small strain kinematics. Sections~\ref{sect-mixedform-small-cont} and~\ref{sect-review-averaging}
present the continuous and discrete problems respectively, while Section
\ref{sect-review-stabilization} closes this part with comments on the main results presented by Lamichhane
\cite{lamichhane2009hu} and the need for a stabilization term. 

\subsection{Mixed variational formulation}
\label{sect-mixedform-small-cont}

\begin{figure}[t]
    \centering
    \includegraphics[width=\textwidth]{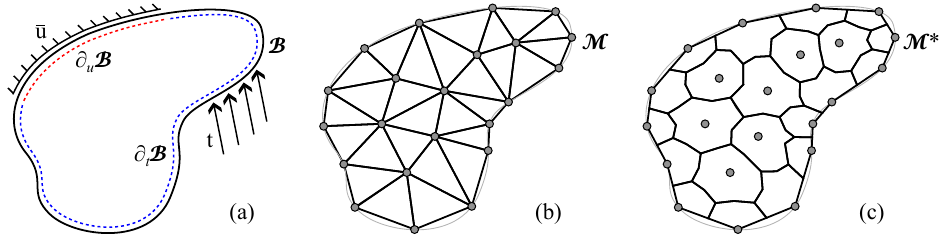}
    \caption{Mixed formulation. (a) Continuum setting. (b) Primal mesh. (c) Dual mesh.}
    \label{fig-cont-meshes}
\end{figure}

Let $\mathcal{B} \subset \mathbb{R}^d$ be a bounded domain occupied by a body of dimension $d=2$ or~$3$,
with boundary $\partial\mathcal{B} = {\partial_u\mathcal{B}} \cup  {\partial_t\mathcal{B}}$ and
${\partial_u\mathcal{B}}\cap {\partial_t\mathcal{B}} = \emptyset$. The displacement of the body is
the vector field $\bm{u} : \mathcal{B} \rightarrow \mathbb{R}^d$. The body is subject to body
forces $\bm{b}:\mathcal{B}\rightarrow \mathbb{R}^d$, surface tractions
$\bm{t}:{\partial_t\mathcal{B}}\rightarrow \mathbb{R}^d$, and imposed displacements $\bm{u} =
\Bar{\bm{u}}$ on ${\partial_u\mathcal{B}}$. This setting is illustrated in
Figure~\ref{fig-cont-meshes}(a). To write the mechanical problem in variational form, let us define
the functional space
\begin{equation}
    \mathcal{U} := \{ \bm{u} \in [H^1(\mathcal{B})]^d \ | \ \bm{u} = \Bar{\bm{u}} \text{ on } {\partial_u\mathcal{B}} \} \ ,
\end{equation}
where $[H^1(\mathcal{B})]^d$ is the Hilbert space of $d$-dimensional vector fields defined
on $\mathcal{B}$ that are square integrable, the space $[L^2(\mathcal{B})]^d$, and whose gradients are also square integrable.
With these ingredients, we recall that the equilibrium displacement field satisfies
\begin{equation}
    \inf_{\bm{u} \in \mathcal{U}}  J({\bm{u}})\,,
    \label{eq-minmechanical}
\end{equation}
with potential energy $J: \mathcal{U}\to \mathbb{R}$ defined as:
\begin{equation}
    J({\bm{u}})= \int_\mathcal{B} {\phi({\nabla^s \bm{u}})} \ dV - \int_\mathcal{B} {\bm{b} \cdot \bm{u}} \ dV - \int_{\partial_t\mathcal{B}}{\bm{t} \cdot \bm{u}} \ dA .
\end{equation}
The first integral in this functional corresponds to the internal potential energy and
it is based on the strain energy function $\phi$ which itself is a function of the symmetric part of the displacement
gradient, i.e., $\nabla^s\bm{u}$. The stationarity condition of this problem reads
\begin{equation}
    0 =
    \int_\mathcal{B} {
    D\phi({\nabla^s \bm{u}}) : \nabla^s \delta\bm{u} } \ dV 
    -
    \int_\mathcal{B} {\bm{b} \cdot \delta\bm{u}} \ dV 
    -
    \int_{\partial_t\mathcal{B}}{\bm{t} \cdot \delta\bm{u}} \ dA \ , 
    \label{eq-bilinearCmch}
\end{equation}
for all $\delta\bm{u} \in \mathcal{U}_0$, the space of admissible variations defined as
\begin{equation}
    \mathcal{U}_0 := \{ \delta\bm{u} \in [H^1(\mathcal{B})]^d \ | \ 
    \delta\bm{u} = \bm{0} \text{ on } {\partial_u\mathcal{B}} \} .
\end{equation}

Problem~\eqref{eq-bilinearCmch} is the starting point of the displacement-based finite element method. 
For the nodal strain formulation presented in this paper, the minimization problem~\eqref{eq-minmechanical} is rewritten as a saddle point problem based on the Hu-Washizu variational
principle~\cite{washizu1982vj}. To introduce this new problem, first we define a function space for the strain and stress fields, namely,
\begin{equation}
    \mathcal{E} := \{ \bm{\varepsilon} \in [L^2(\mathcal{B})]^{d\times d} \ | \ \bm{\varepsilon} \text{ is symmetric } \}
    .
\end{equation}
Then, the saddle point problem is 
\begin{equation}
\label{eq-mixed-L}
    \inf_{\bm{u} \in \mathcal{U}, \ \bm{\varepsilon} \in \mathcal{E}
    } \sup_{\bm{\sigma} \in \mathcal{E}} L({\bm{u}},{\bm{\varepsilon}},{\bm{\sigma}}) , 
\end{equation}
with the functional $L:\mathcal{U}\times \mathcal{E} \times \mathcal{E}\to \mathbb{R}$ defined as
\begin{equation}
  L({\bm{u}},{\bm{\varepsilon}},{\bm{\sigma}}) :=
  \int_\mathcal{B} {\phi({\bm{\varepsilon}})} \ dV + \int_\mathcal{B} {\bm{\sigma} : (\nabla^s \bm{u}-\bm{\varepsilon})} \ dV + \alpha \int_\mathcal{B} {\hat{\phi}({\nabla^s \bm{u}-\bm{\varepsilon}})} \ dV - \int_\mathcal{B} {\bm{b} \cdot \bm{u}} \ dV - \int_{\partial_t\mathcal{B}}{\bm{t} \cdot \bm{u}} \ dA \ .
    \label{eq-HWenergy} 
\end{equation}

In this expression, we have introduced a positive constant $\alpha$ and an ancillary stored energy $\hat{\phi}$ 
defined on symmetric, second-order tensors $\bm{d}$, and of the form
\begin{equation}
    \hat{\phi}({\bm{d}}) := \frac{1}{2} \bm{d}: \hat{\mathbb{C}}\bm{d} \ ,
    \label{eq-stabss}
\end{equation}
with $\hat{\mathbb{C}}$ being an auxiliary stiffness tensor $\hat{\mathbb{C}} := \hat{\lambda} \bm{1}\otimes\bm{1} + 2\hat{\mu}\mathbb{I}$, where 
$(\hat{\lambda},\hat{\mu})$ play the role of Lam\'e parameters. The stationarity conditions of problem~\eqref{eq-mixed-L} are:
find $(\bm{u},\bm{\varepsilon},\bm{\sigma}) \in \mathcal{U}\times\mathcal{E}\times\mathcal{E}$ such that
\begin{subequations}
    \begin{align}
        \int_\mathcal{B} {\bm{\sigma}:\nabla^s \delta\bm{u}} \ dV + \alpha \int_\mathcal{B} {\hat{\mathbb{C}} (\nabla^s \bm{u}-\bm{\varepsilon}):\nabla^s \delta\bm{u}} \ dV
        - \int_\mathcal{B} {\bm{b} \cdot \delta\bm{u}} \ dV - \int_{\partial_t\mathcal{B}}{\bm{t} \cdot \delta\bm{u}} \ dA &= 0,
        \label{eq-hw-1}
        \\
        \int_\mathcal{B} {\left(D \phi({\bm{\varepsilon}}) - \bm{\sigma}\right): \delta\bm{\varepsilon}} \ dV - \alpha \int_\mathcal{B} {\hat{\mathbb{C}} (\nabla^s \bm{u}-\bm{\varepsilon}):\delta\bm{\varepsilon}} \ dV &= 0,
        \label{eq-hw-2}
        \\
        \int_\mathcal{B} {(\nabla^s \bm{u} - \bm{\varepsilon}):\delta\bm{\sigma}} \ dV &= 0,
        \label{eq-hw-3}
    \end{align}%
    \label{eq-huwashizu}%
\end{subequations}%
for all $(\delta\bm{u},\delta\bm{\varepsilon},\delta\bm{\sigma}) \in \mathcal{U}_0\times\mathcal{E}\times\mathcal{E}$. 
Since Eq.~\eqref{eq-hw-3} weakly imposes $\bm{\varepsilon}=\nabla^s\bm{u}$, the integrals 
proportional to~$\alpha$ in Eqs.~\eqref{eq-hw-1} and ~\eqref{eq-hw-2} identically 
vanish, showing that the auxiliary energy plays no role whatsoever in the problem. Later, we will show that this is not the case for the discrete problem, where
this term actually can be selected to stabilize a certain finite element formulation.

\subsection{Finite element discretization}
\label{sect-review-averaging}

Consider a conforming partition of the domain $\mathcal{B}$ into simplicial elements or
quad/hex-type elements. This partition defines a primal mesh 
$\mathcal{M} = \{ e_\alpha \}$ of cells $e_\alpha$, with vertex set
$\mathcal{N} = \{ n_\beta \}$ of cardinality $N_n$, and a dual mesh
${\mathcal{M}^*} = \{ e^\beta \}$ where each cell $e^\beta$ is associated with node $n_\beta$. Details
of how to construct the dual mesh can be found in the literature \cite{dohrmann2000node, lamichhane2009hu}. The
relation between the primal and dual meshes is illustrated in Figure~\ref{fig-cont-meshes}(b)
and~(c), respectively. Based on these partitions of the body, we can define the function spaces
\begin{equation}
    \begin{split}
      \mathcal{U}^h &:= \{ \bm{u}^h \in \mathcal{U},\  \bm{u}^h|_{e_\alpha} \in [R_1(e_\alpha)]^d ,
                      \;e_\alpha \in \mathcal{M} \} \ , \\
        \mathcal{U}_0^h &:=\{ \bm{u}^h \in \mathcal{U}_0, \  \bm{u}^h|_{e_\alpha} \in [R_1(e_\alpha)]^d , \;e_\alpha \in \mathcal{M} \} \ , \\
        \mathcal{E}^h &:=\{ \bm{\varepsilon}^h \in \mathcal{E}, \ \bm{\varepsilon}^h|_{e^\beta} \in [R_0(e^\beta)]^{d\times d} ,\; e^\beta \in {\mathcal{M}^*} \} \ ,
    \end{split}
    \label{eq-spacesFE0}
\end{equation}
where $R_m(e_{\alpha})$ is $P_m(e_{\alpha})$, the space of polynomials of degree $m$ if $e_{\alpha}$
is a simplicial element, or $Q_m(\alpha)$, the space of isoparametric maps of tensor-product
functions if $e_{\alpha}$ is a quad/hexahedron. Then, the finite element approximation to problem~\eqref{eq-huwashizu} is as follows: find $(\bm{u}^h,\bm{\varepsilon}^h,\bm{\sigma}^h) \in \mathcal{U}^h\times\mathcal{E}^h\times\mathcal{E}^h$ such that
\begin{subequations}
    \begin{align} 
        \int_\mathcal{B} {\bm{\sigma}^h:\nabla^s \delta\bm{u}^h} \ dV + \alpha \int_\mathcal{B} {\hat{\mathbb{C}} (\nabla^s \bm{u}^h-\bm{\varepsilon}^h):\nabla^s \delta\bm{u}^h} \ dV
        - \int_\mathcal{B} {\bm{b} \cdot \delta\bm{u}^h} \ dV - \int_{\partial_t\mathcal{B}}{\bm{t} \cdot \delta\bm{u}^h} \ dA &= 0,  \label{eq-huwashizu-discr1}\\
        \int_\mathcal{B} {\left( D \phi({\bm{\varepsilon}^h}) - \bm{\sigma}^h\right): \delta\bm{\varepsilon}^h} \ dV - \alpha \int_\mathcal{B} {\hat{\mathbb{C}} (\nabla^s \bm{u}^h-\bm{\varepsilon}^h):\delta\bm{\varepsilon}^h} \ dV &= 0, \label{eq-huwashizu-discr2}\\
        \int_\mathcal{B} {(\nabla^s \bm{u}^h - \bm{\varepsilon}^h):\delta\bm{\sigma}^h} \ dV &= 0, \label{eq-huwashizu-discr3}
    \end{align}\label{eq-huwashizu-discr}%
\end{subequations}%
for all $(\delta\bm{u}^h,\delta\bm{\varepsilon}^h,\delta\bm{\sigma}^h) \in
\mathcal{U}_0^h\times\mathcal{E}^h\times\mathcal{E}^h$.

\begin{figure}[ht]
    \centering
    \includegraphics[width=12.5cm]{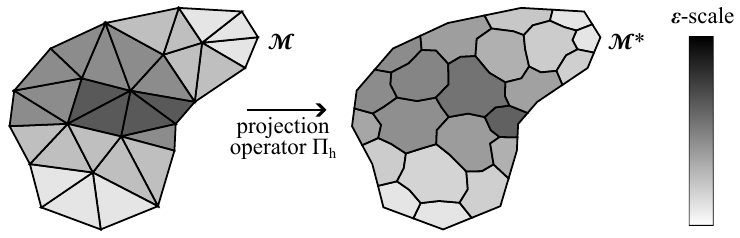}
    \caption{Action of the projection operator $\Pi_h:\mathcal{E} \rightarrow \mathcal{E}^h$.}
    \label{fig-projection}
\end{figure}

Next, let us simplify Eq.~\eqref{eq-huwashizu-discr3} using the properties of the interpolation
spaces. For that, let $\chi^\beta$ be the characteristic function of the dual cell $e^\beta$, that is, the function
\begin{equation}
    \chi^\beta (\bm{x}) :=
    \begin{cases}
        1 \quad \text{if} \ \bm{x} \in e^\beta \ , \\
        0 \quad \text{otherwise} \ .
    \end{cases}
\end{equation}
Then, choose $\delta\bm{\sigma}^h = \chi^\beta \ \bm c$ in Eq.~\eqref{eq-huwashizu-discr3}, 
where $\bm c$ is an arbitrary constant tensor. From the resulting relation it is possible to introduce the projection operator $\Pi_h:\mathcal{E} \rightarrow \mathcal{E}^h$ that links $\nabla^s\bm{u}^h$ with $\bm{\varepsilon}^h$, as in
\begin{equation}
    \bm{\varepsilon}^h := \Pi_h(\nabla^s\bm{u}^h) := \sum_{\beta=1}^{N_n} \bm{\varepsilon}^h_\beta \ \chi^\beta \ ,
    \label{eq-projectionOper}
\end{equation}
where
\begin{equation}
    \bm{\varepsilon}^h_\beta := \bm{\varepsilon}^h|_{e^\beta} := \frac{1}{|e^\beta|} \int_{e^\beta} {\nabla^s \bm{u}^h} \ dV \ .
    \label{eq-dualstrain}
\end{equation}
We note that $\bm{\varepsilon}^h_\beta$ is constant on each cell of the dual mesh. The projection operation is illustrated in Figure~\ref{fig-projection} for a triangular mesh.

For economy of notation, in the following, the projection $\Pi_h(\nabla^s\bm{u}^h)$ is replaced by
${\nablabar^s}\bm{u}^h$, which reflects the fact that this object is computed from the
displacement field using a modification of the strain operator. The range of $\nablabar^s$ is the
space of symmetric tensors that are constant on each dual cell. This is a closed space and thus, by
the orthogonal decomposition theorem applied to $\mathcal{E}^h$, for all $\bm{v}^h\in [H^1(\mathcal{B})]^d$,
we have that
\begin{equation}
\begin{split}
    \int_\mathcal{B} \nablabar^s\bm{u}^h:\nabla^s \bm{v}^h \ dV &= \int_\mathcal{B} {\left(\sum_{\beta=1}^{N_n} \bm{\varepsilon}^h_\beta(\bm{u}^h) \ \chi^\beta\right):\nabla^s \bm{v}^h} \ dV 
    \\
    &= \sum_{\beta=1}^{N_n} \bm{\varepsilon}^h_\beta(\bm{u}^h) : \int_{e^\beta} {\nabla^s \bm{v}^h} \ dV 
    \\
    &= \sum_{\beta=1}^{N_n} \bm{\varepsilon}^h_\beta(\bm{u}^h) : \bm{\varepsilon}^h_\beta(\bm{v}^h) |e^\beta| \\
    &= \int_\mathcal{B} {{\nablabar^s}\bm{u}^h:{\nablabar^s}\bm{v}^h} \ dV .
\end{split}
\label{eq-projOpProp}
\end{equation}
Using this property, it is possible to eliminate the second term from Eq.~\eqref{eq-huwashizu-discr2}.

Repeating the same steps employed in the definition of the assumed strain, one obtains a representation for the stress field analogous to that of the strains as
\begin{equation}
    \bm{\sigma}^h = \sum_{\beta=1}^{N_n} \bm{\sigma}^h_\beta \ \chi^\beta \ ,
    \label{eq-stressfield}
\end{equation}
where
\begin{equation}
    \bm{\sigma}^h_\beta := \bm{\sigma}^h|_{e^\beta} := \frac{1}{|e^\beta|} \int_{e^\beta} {D \phi({\bm{\varepsilon}^h})} \ dV \ .
\end{equation}
This identity shows that the fictitious stored energy $\hat{\phi}$ does not affect the constitutive relation for the assumed stress field, which is computed exclusively from the free energy $\phi({\bm{\varepsilon}})$. Remarkably, if the cell ${e^\beta}$ is homogeneous, the definition of the algorithmic stress simplifies to
\begin{equation}
    \bm{\sigma}^h_\beta = D \phi({\bm{\varepsilon}^h_\beta}) \ .
    \label{eq-stress}
\end{equation}

Collecting the result~\eqref{eq-projOpProp} and the representations~\eqref{eq-projectionOper}
and~\eqref{eq-stressfield}, problem~\eqref{eq-huwashizu-discr} can be stated alternatively as follows: find $\bm{u}^h \in \mathcal{U}^h$ such that
\begin{equation}
\begin{split}
    \int_\mathcal{B} {\bm{\sigma}^h(\bm{u}^h):{\nablabar^s}\delta\bm{u}^h} \ dV - \alpha \int_\mathcal{B} {\hat{\mathbb{C}} \ {\nablabar^s}\bm{u}^h:{\nablabar^s}\delta\bm{u}^h} \ dV + \alpha \int_\mathcal{B} {\hat{\mathbb{C}} \ \nabla^s \bm{u}^h:\nabla^s \delta\bm{u}^h} \ dV = \\
    = \int_\mathcal{B} {\bm{b} \cdot \delta\bm{u}^h} \ dV + \int_{\partial_t\mathcal{B}}{\bm{t} \cdot \delta\bm{u}^h} \ dA , 
\end{split}
\label{eq-nodalFEM}
\end{equation}
for all $\delta\bm{u}^h \in \mathcal{U}_0^h$.

\begin{remark}
The only unknown in problem~\eqref{eq-nodalFEM} is the displacement field $\bm{u}^h \in \mathcal{U}^h$. In
other words, although the formulation starts from a mixed problem, the choice of function spaces for
the strain and stress fields leads to a reduced problem stated entirely in terms of displacements
that are piecewise polynomials on the primal mesh, just like a standard displacement-based finite
element discretization. This property is shared by other \emph{assumed strain} methods
\cite{wilson1973ty,simo1990uj,bischoff2007az,romero-bubbles-2007}. 
\end{remark}

\begin{remark}
  The integrals involving the averaged strain $\nablabar^s$ are computed over the dual mesh, whereas those
  involving only the displacement field and the consistent gradient $\nabla^s$ are computed over the primal mesh. Thus, the bilinear form
  of the nodal strain method, given by the left-hand side of Eq.~\eqref{eq-nodalFEM}, contains two
  terms evaluated on the dual mesh and one term evaluated on the primal mesh. If the stabilization
  terms are ignored, problem~\eqref{eq-nodalFEM} has the same structure as
  problem~\eqref{eq-bilinearCmch}, except that the discrete bilinear form differs from its
  continuous counterpart through the use of the projection operator.
\end{remark}

\begin{remark}
  In connection with the previous remark, some misconceptions exist in the literature regarding the
  continuity of the strain and stress fields. Some works claim that the strain and stress fields
  admit a continuous representation obtained by interpolating their \emph{nodal values} with the
  shape functions of the primal mesh~\cite{dohrmann2000node,castellazzi2015linear}. Although such an
  interpolation may be acceptable as a post-processing step, the definition of the functional spaces
  in Eq.~\eqref{eq-spacesFE0} shows that the strain and stress fields are approximated as piecewise
  constant quantities over the dual partition of the domain. In this sense, rather than \emph{nodal
    quantities}, the term \emph{dual-cell quantities} may be more appropriate.
\end{remark}

\begin{remark}
  The distinction between \emph{nodal} and \emph{cell} quantities turns out to be important for the
  calculation and numerical treatment of history variables of inelastic materials. This is because
  computations involving inelastic materials often make use of \emph{history variables} that encode,
  at each instant of the analysis, the effect of the past deformation (see, e.g.,
  \cite{coleman1967di,simo2000vm,romero2026eh}). In numerical discretizations, this information is
  associated to material points. In standard finite element discretizations, for example, material
  points are placed at the quadrature points, and this is where the strains and stresses are
  evaluated. In nodal finite element methods, material points correspond to dual cells
  and can be bijectively associated with the nodes.
  \end{remark}

Regardless of terminology, the key point is that in the formulation~\eqref{eq-huwashizu-discr}, the
correspondence between dual cells and nodes makes it possible to associate strain, stress, and
material history with the mesh nodes. From this perspective, the resulting nodal structure of the
method admits a twofold interpretation. On the one hand, it can be understood as a finite element
formulation based on properly defined function spaces, with piecewise constant strain and stress
fields associated with nodal patches. On the other hand, it can be interpreted as a strict
particle method in which the complete material information is stored at, or naturally
associated with, the nodes, which therefore behave as material points. In this sense, even though
a material model must also be defined on the primal mesh because of the stabilization term in
Eq.~\eqref{eq-nodalFEM}, this auxiliary material does not affect the actual constitutive response
in Eq.~\eqref{eq-stress}. Moreover, since this ancillary material is elastic, it does not
introduce any additional material history.

\subsection{On the need for stabilization and Lamichhane's results}
\label{sect-review-stabilization}

Lamichhane's work \cite{lamichhane2009hu} contributed significantly to the theoretical foundations
of the nodal strain finite element method. For linear elasticity, with the real material defined by
the parameters $(\lambda,\mu)$ and using tetrahedral elements, he proved the coercivity of the bilinear form
in~\eqref{eq-nodalFEM} when choosing $(\hat{\lambda},\hat{\mu})=(0,\mu)$ for the stabilizing material. He also established
optimal convergence of order $\mathcal{O}(h)$ in the $H^1$-norm for the displacement field. Furthermore, he showed
that the coercivity constant is independent of the ratio $\lambda/\mu$ of the Lamé parameters, which
explains the locking-free behaviour of the method, and discussed the choice $\alpha=1$ as a natural
compromise between preserving coercivity and avoiding pollution of the approximate solution.

Other works have also addressed the choice of material parameters in the stabilization term
\cite{puso2006stabilized, gee2009uniform}. In particular, the isochoric stabilization proposed in
\cite{gee2009uniform} is, in the linear case, equivalent to choosing $\hat{\mu}=\mu$ and
$\hat{\kappa}=\hat{\lambda}+2\hat{\mu}/3=0$. Alternative stabilization procedures have also been proposed. For instance, in the
context of NSFEM, \cite{tong2018high, chen2019gradient} introduced a mechanism based on the
assumption of a continuous and differentiable strain field. From our point of view, this introduces
a variational inconsistency, since it is not compatible with the function spaces defined previously.

Despite the efforts devoted to obtaining a stable numerical procedure, the loss of coercivity in the
absence of a stabilization term still appears to be insufficiently recognized. Some works have
usefully illustrated possible zero-energy deformation modes in this family of methods, particularly
in the context of infinite domains \cite{puso2008meshfree, sivapuram2019energy}. However, it is
often argued that boundaries provide a certain \emph{stabilizing effect}, so that meshes with high
nodal-surface-to-volume ratios, i.e. coarser meshes, appear more stable \cite{puso2008meshfree,
  artioli2014assumed}. It has also been suggested that instabilities are mainly relevant in dynamic
analyses \cite{zhang2010temporal, feng2014temporal, sivapuram2019energy}.

This section closes with a simple example which shows that the stabilization term is necessary even in a stationary
analysis on finite domains with coarse meshes.

\subsubsection{Numerical example}
\label{subs-stability-small-example}

\begin{figure}[t]
    \centering
    \includegraphics[width=6cm]{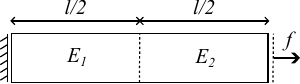}
    \caption{Problem setting.}
    \label{fig-stabscheme}
\end{figure}

The example under consideration consists of a straight bar of length $l$ and square cross section of
area $S$, with one end fixed and an axial load of magnitude $f$ applied at the other end. Each half of the bar is composed of an elastic material
characterized by Young's moduli $E_1$ and $E_2$ and Poisson's ratio $\nu$ (see Figure~\ref{fig-stabscheme} for an
illustration). The load is applied uniformly over the cross section. To obtain a purely one-dimensional solution, we set Poisson's ratio to $\nu = 0$. Then,
the analytical solution of the displacement field in the $x$-direction --- the only relevant one for the one-dimensional solution --- is
\begin{equation}
    \Tilde{u}_x = 
    \begin{cases}
        \hat{\varepsilon}_1\,x + \hat{\varepsilon}_1 \ l/2 \quad &\text{if} \ \ -l/2\leq x \leq 0\\
        \hat{\varepsilon}_2\,x + \hat{\varepsilon}_1 \ l/2 \quad &\text{if} \ \ 0 < x\leq l/2
    \end{cases}
\end{equation}
with $\hat{\varepsilon}_i := f/(SE_i)$, and $i=1,2$, being dimensionless constants. In addition, $\Tilde{u}_y = \Tilde{u}_z = 0$. For this simple problem, the
strain component $\Tilde{\varepsilon}_{xx}$
is, therefore, piecewise constant:
\begin{equation}
    \Tilde{\varepsilon}_{xx} = 
    \begin{cases}
        \hat{\varepsilon}_1 \quad &\text{if} \ \ -l/2\leq x \leq 0 \\
        \hat{\varepsilon}_2 \quad &\text{if} \ \ 0 < x\leq l/2
    \end{cases}
\end{equation}

\begin{figure}[htbp]
    \centering

    \begin{subfigure}{0.32\textwidth}
        \centering
        \includegraphics[width=\linewidth]{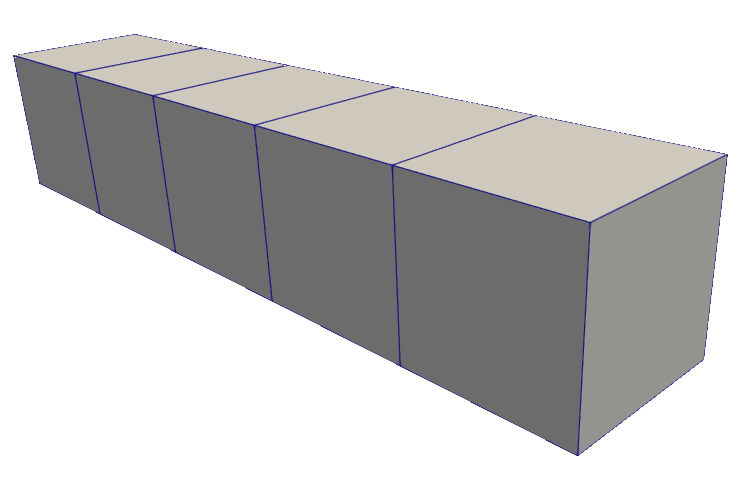}
    \end{subfigure}
    \hfill
    \begin{subfigure}{0.32\textwidth}
        \centering
        \includegraphics[width=\linewidth]{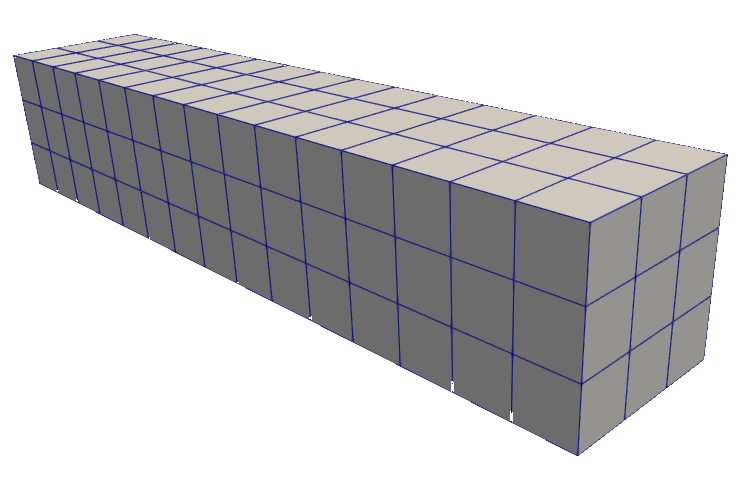}
    \end{subfigure}
    \hfill
    \begin{subfigure}{0.32\textwidth}
        \centering
        \includegraphics[width=\linewidth]{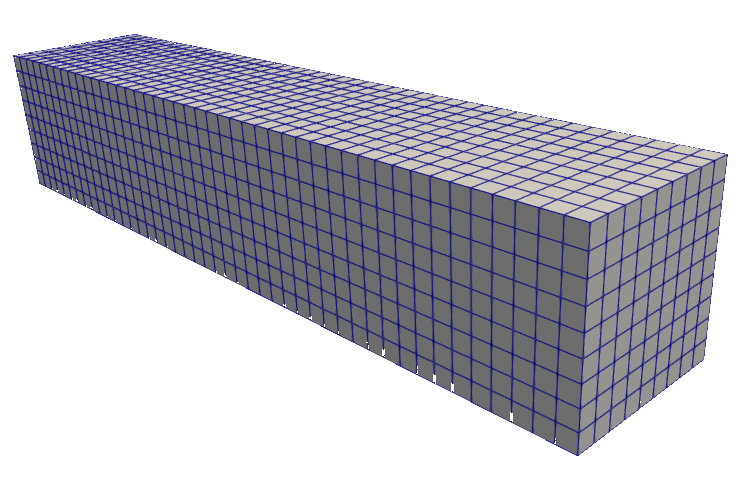}
    \end{subfigure}
    \caption{Sequence of meshes used in the two-material bar example.}
    \label{fig-stabexample-meshes}
\end{figure}

To study this problem numerically, we select $f=0.1$, $l=5$, $S=1$, $E_1=1$, $E_2=2$, and then proceed to obtain
numerical approximations with three meshes of increasingly smaller elements. We employ the nodally
integrated method with and without stabilization, the latter corresponding to $\alpha=0$ in problem~\eqref{eq-nodalFEM}. For the stabilized case we define a single stabilizing material for the whole body, with properties following Lamichhane's approach: $\hat{\lambda} = 0$, $\hat{\mu} = \mu_1 = E_1/(2(1+\nu))=0.5$, and $\alpha=1.0$. The first mesh under
consideration consists of $5\times1\times1$ hexahedral elements, and each subsequent mesh is obtained by subdividing each
element into $3\times3\times3$ new ones. This sequence of meshes is depicted in Figure~\ref{fig-stabexample-meshes}.

The sequence of meshes is designed in this way because of the manner in which the method represents a heterogeneous body. Recall that Eq.~\eqref{eq-stress} is based on the assumption that the cells of the dual mesh are homogeneous, which implies that the interface between materials lies inside the elements of the primal mesh. Interestingly, this emphasizes the interpretation of the approach as a particle method: the constitutive information is stored at the nodes (or dual cells). Therefore, nodes cannot lie on the interface between materials. This, in turn, prevents the method from reproducing the exact solution for this problem: the exact displacement field $\Tilde{\bm{u}}$ (approximated on the primal mesh) changes its description from one material to the other, which cannot be captured by a single element shape function.

\begin{figure}[htbp]
    \centering

    \begin{subfigure}{\textwidth}
        \centering
        \begin{minipage}{0.48\textwidth}
            \centering
            \includegraphics[width=\linewidth]{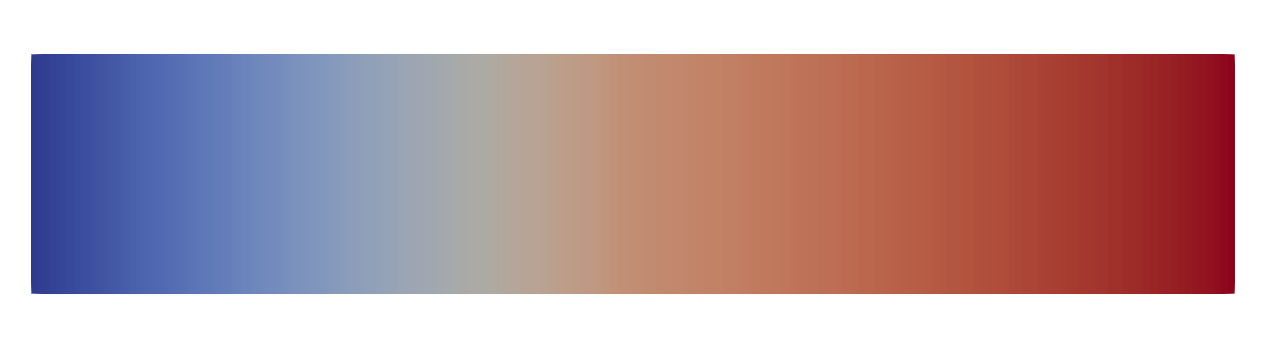}
        \end{minipage}
        \hfill
        \begin{minipage}{0.48\textwidth}
            \centering
            \includegraphics[width=\linewidth]{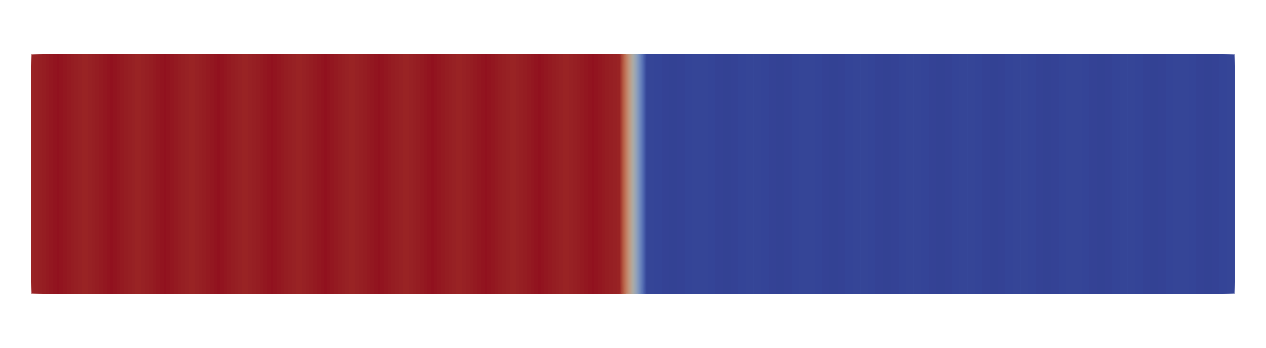}
        \end{minipage}
        \caption{without stabilization}
        \label{fig-stabexample-without}
    \end{subfigure}

    \begin{subfigure}{\textwidth}
        \centering
        \begin{minipage}{0.48\textwidth}
            \centering
            \includegraphics[width=\linewidth]{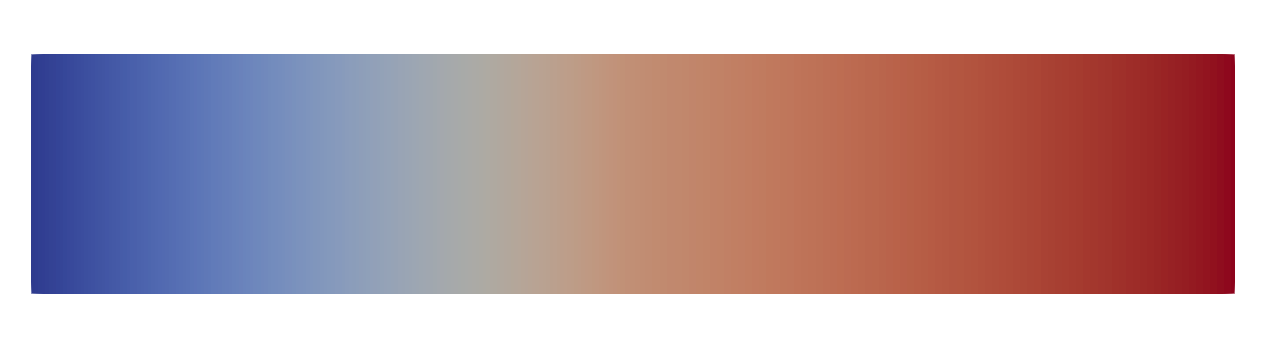}
        \end{minipage}
        \hfill
        \begin{minipage}{0.48\textwidth}
            \centering
            \includegraphics[width=\linewidth]{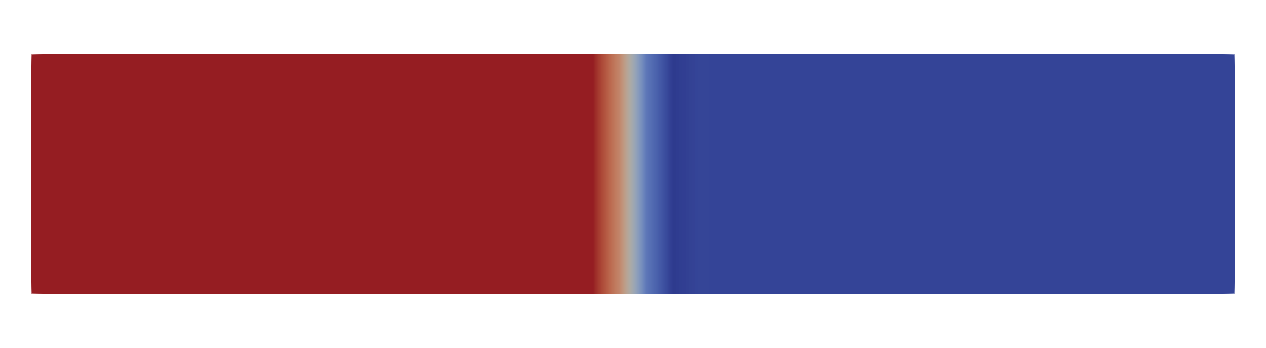}
        \end{minipage}
        \caption{with stabilization}
        \label{fig-stabexample-with}
    \end{subfigure}

    \begin{minipage}{0.48\textwidth}
        \centering
        \includegraphics[width=\linewidth]{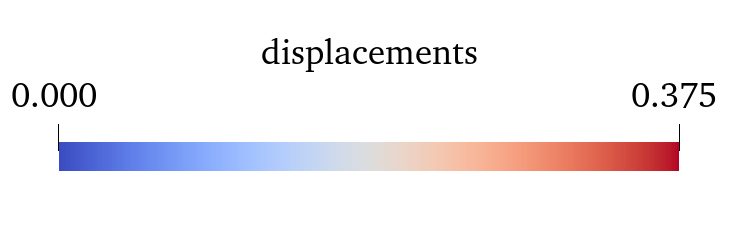}
    \end{minipage}
    \hfill
    \begin{minipage}{0.48\textwidth}
        \centering
        \includegraphics[width=\linewidth]{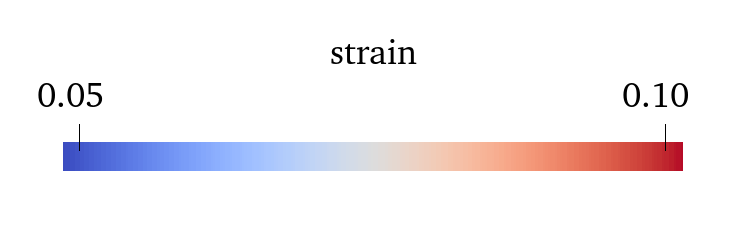}
    \end{minipage}

    \caption{Results obtained by nodal strain finite element method (a) without, and (b) with stabilization, for the $45\times 9 \times 9$ mesh.}
    \label{fig-stabexample}
\end{figure}

As an example, Figure~\ref{fig-stabexample} shows the results obtained for displacements ${u}_x$ and strains ${\varepsilon}_{xx}$ for the finest mesh. In line with the discussion at the end of
Section~\ref{sect-review-averaging}, the strain field (right column of plots) is actually piecewise
constant, with values that can be linked to the nodes of the mesh. The displayed results show a
linearly interpolated field for convenience in the representation. The solution generated by the non-stabilized method exhibits oscillations in both displacements and strains, although these are only clearly visible in the latter.

\begin{figure}[htbp]
    \centering
    \includegraphics[width=15cm]{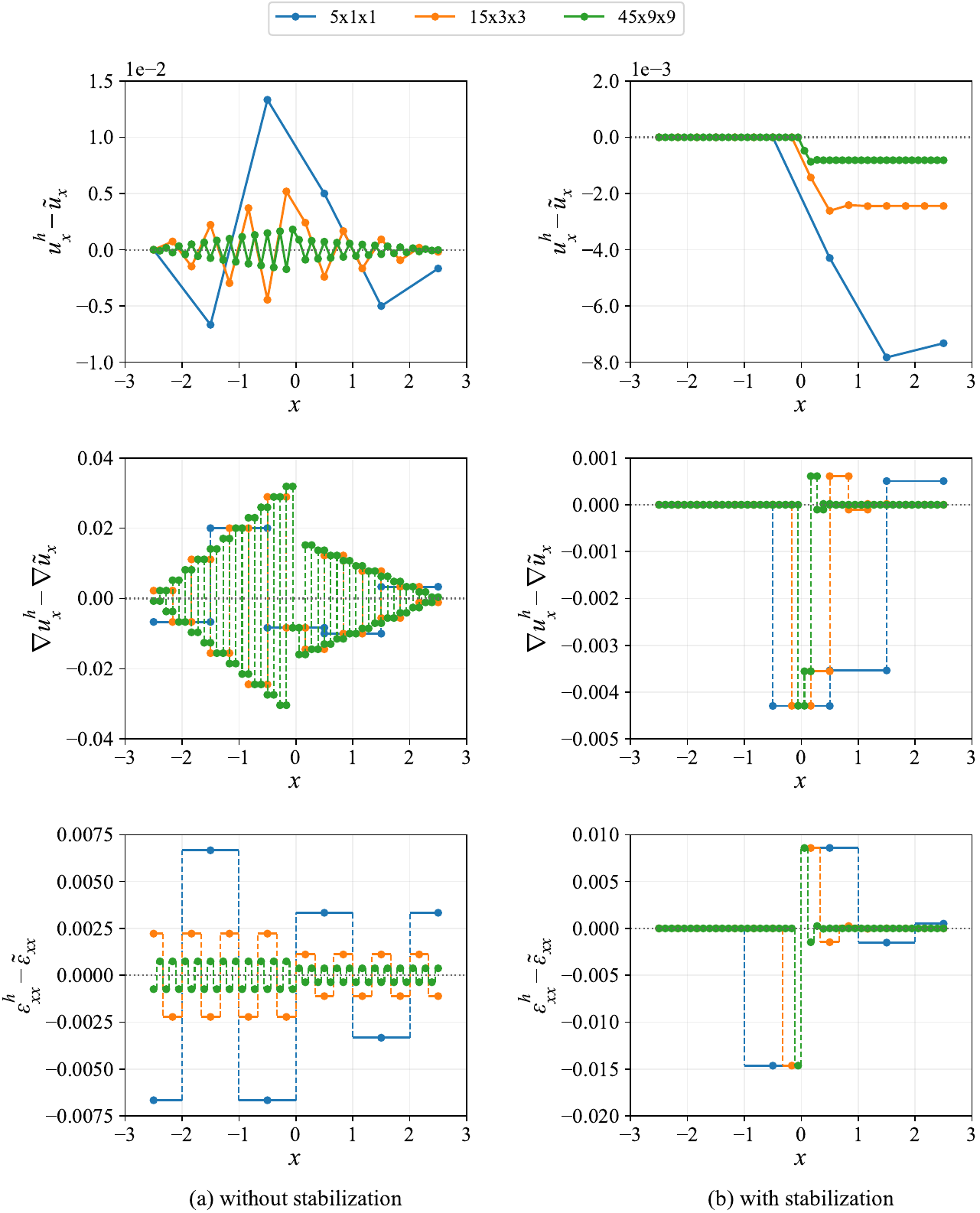}
    \caption{Errors obtained for the nodal strain FEM: (a) without stabilization, (b) with stabilization. Top: errors in displacements. Centre: error in displacement gradients. Bottom: error in strains (dual mesh). Circles indicate nodes of the mesh.}
    \label{fig-results-stab}
\end{figure}

Further insight into the behaviour of the numerical method can be obtained from Figure~\ref{fig-results-stab}. This figure presents the errors in displacements, displacement-compatible gradients, and strains along the axis of the bar.
It can be observed that, irrespective of the mesh size, the absence of the stabilization term
produces oscillations in both the displacement and strain fields. Interestingly, despite these oscillations, both displacements and strains converge in the $L^2$-norm. However, the displacement gradients do not converge, which means that the displacement field does not converge in the $H^1$-norm. As expected, convergence is recovered when stabilization is included.

This example shows that, for this case, the method does not converge in the expected norm in the absence of stabilization. Furthermore, unrealistic oscillations appear in the context of a static analysis with coarse meshes.

\section{Extension to large-strain kinematics}
\label{sect-mixedform-large}
This section presents the main theoretical contribution of the paper. The Hu–Washizu
functional~\eqref{eq-HWenergy} is generalized to the finite strain range, allowing the formulation
to accommodate large deformations. The continuum and discrete problems are established in
Sections~\ref{sect-mixedformprob-large} to~\ref{sect-mixedformprob-largeFE}, including a generalization of the stabilization term in Section~\ref{sect-mixedformprob-large-stab}.
Section~\ref{sect-projop-reference} introduces the generalization of the assumed strain operator,
accounting for changes of the reference configuration. Section~\ref{sect-stablarge} closes with a discussion
of the choice of the stabilizing function in the large-strain scenario.

\subsection{Continuum formulation}
\label{sect-mixedformprob-large}

Let ${\mathcal{B}_0} \subset \mathbb{R}^d$ be a bounded domain occupied by the reference configuration of a deformable body, with
boundary $\partial{\mathcal{B}_0} = {\partial_u{\mathcal{B}_0}} \cup {\partial_t{\mathcal{B}_0}}$ and
${\partial_u{\mathcal{B}_0}}\cap {\partial_t{\mathcal{B}_0}} = \emptyset$. In addition, let
$\mathcal{B} \subset \mathbb{R}^d$ be the current configuration of the same body which is related to the reference configuration
by the deformation $\bm{\varphi} : {\mathcal{B}_0} \rightarrow \mathcal{B}$ with deformation gradient
$\bm{F}=\nabla_0\bm{\varphi}$, where $\nabla_0$ denotes the material gradient
operator. The body is subject to body forces $\bm{b}_0:{\mathcal{B}_0}\rightarrow \mathbb{R}^d$ per unit reference volume, surface tractions
$\bm{t}_0:{\partial_t{\mathcal{B}_0}}\rightarrow \mathbb{R}^d$ in the reference configuration, and imposed deformation
$\bm{\varphi} = \Bar{\bm{\varphi}}$ on ${\partial_u{\mathcal{B}_0}}$. As in the small strain setting, we strive to formulate the mechanical equilibrium
variationally and, for that, we introduce the appropriate functional spaces 
for the deformation, the deformation gradient, and the Piola-Kirchhoff stress. 
From the analysis point of view, this framework is considerably more complex than the linear case (see, e.g., \cite{ciarlet1988ux}). Hence, we introduce the following functional spaces that are 
commonly used in finite elasticity, without claiming that these
are enough to prove the existence of a solution to the problem:
\begin{equation}
    \begin{split}
        \mathcal{X} :=& \ \{ \bm{\varphi} \in [W^{1,p}({\mathcal{B}_0})]^d \ | \ \bm{\varphi} = \Bar{\bm{\varphi}} \text{ on } {\partial_u\mathcal{B}}_0 \} \ , \\
        \mathcal{X}_0 :=& \ \{ \bm{\varphi} \in [W^{1,p}({\mathcal{B}_0})]^d \ | \ \bm{\varphi} = \bm{0} \text{ on } {\partial_u\mathcal{B}}_0 \} \ , \\
        \mathcal{F} :=& \ [L^2({\mathcal{B}_0})]^{d\times d} \ .
    \end{split}
\end{equation}
Here, $W^{1,p}(\mathcal{B}_0)$, with $p>1$, is the Sobolev space of functions defined on the reference
configuration that belong to $L^p(\mathcal{B}_0)$ and their gradients too. 
With these ingredients, we can introduce the Hu-Washizu functional $L:\mathcal{X}\times \mathcal{F}
\times \mathcal{F}\to \mathbb{R}$ with expression
\begin{equation}
    \label{eq-HWenergy-large} 
\begin{split}
    L({\bm{\varphi}},{\bm{F}},{\bm{P}}) :=& \int_{\mathcal{B}_0} {\phi({\bm{F}})} \ dV + \int_{\mathcal{B}_0} {\bm{P} : (\nabla_0 \bm{\varphi}-\bm{F})} \ dV + 
    \alpha \int_{\mathcal{B}_0} {\hat{\phi}({\nabla_0 \bm{\varphi},\bm{F}})} \ dV 
    \\
    &- \int_{\mathcal{B}_0} {\bm{b}_0 \cdot \bm{\varphi}} \ dV - \int_{{\partial_t\mathcal{B}}_0}{\bm{t}_0 \cdot \bm{\varphi}} \ dA \ .
    \end{split}
\end{equation}
For simplicity of exposition, we have assumed that the body is homogeneous. The tensor $\bm{P}$ has
dimensions of stress and we will prove later that it can be identified with the first
Piola-Kirchhoff stress tensor.

\subsection{The stabilization function}
\label{sect-mixedformprob-large-stab}

As in the small strain case (see Eq.~\eqref{eq-HWenergy}), the third
term of the functional~\eqref{eq-HWenergy-large} is a stabilization 
integral that does not contribute to the continuum
equations, but will play an important role in the discrete setting. In
contrast with the small-strain case, where $\hat{\phi}$ was chosen as a quadratic form of the difference
$\nabla^s \bm{u}-\bm{\varepsilon}$, in the finite strain setting the stabilizing field needs to be
defined more carefully. To generate a variationally consistent formulation with a
convenient implementation, the term $\hat{\phi}$ must satisfy three requirements:
\begin{itemize}
\item[\textbf{R1:}] It must be an objective function.
\item[\textbf{R2:}] Its first variation with respect to $\nabla_0 \bm{\varphi}$ must be separable in the
  discrete setting, after using the properties of the projection operator, into one integral
  evaluated on the dual mesh and one integral evaluated on the primal mesh.
\item[\textbf{R3:}] Its first variation with respect to $\bm{F}$ must not pollute the stress constitutive
  relationship.
\end{itemize}

It is shown next that there is a systematic way to formulate stabilization potentials that satisfy
these three requirements. Let $\hat{g} : \mathbb{R}^{d\times d} \to \mathbb{R}$ be an arbitrary, smooth, frame-indifferent function. The stabilization term
in Eq.~\eqref{eq-HWenergy-large} may be stated as the \emph{Bregman divergence} of $\hat{g}$ \cite{bregman1967relaxation, banerjee2005clustering},
denoted by $D_{B}$ and defined as
\begin{equation}
\hat{\phi}({\nabla_0 \bm{\varphi},\bm{F}})
 \;:=\; D_{B}(\hat{g},\nabla_0 \bm{\varphi},\bm{F})
 \;:=\; \hat{g}(\nabla_0 \bm{\varphi}) - \hat{g}(\bm{F})
        - D\hat{g}(\bm{F}) : (\nabla_0 \bm{\varphi} - \bm{F})\,.
\label{eq-bregman}
\end{equation}
We note that this divergence is just the difference between $\hat{g}$ evaluated at
$\nabla_0 \bm{\varphi}$ and its linearization around $\bm{F}$.

Regarding requirement R1, recall that frame indifference demands the invariance of the
stabilizing energy under a superposed rigid-body motion of the current configuration,
$\bm{\varphi}\mapsto\bm{\varphi}^\ast = \bm{Q} \bm{\varphi} + \bm{c}$, with $\bm{Q} \in SO(d)$ and
$\bm{c}\in \mathbb{R}^d$. Under such a motion, the deformation gradient transforms as
$\nabla_0 \bm{\varphi}^\ast = \bm{Q}\,\nabla_0 \bm{\varphi}$ and, since the independent field
$\bm{F}$ represents the same kinematic object, it is postulated that it transforms in the same way,
$\bm{F}^\ast = \bm{Q} \bm{F}$. Objectivity of Eq.~\eqref{eq-bregman} is then inherited from that of
$\hat{g}$: differentiating the identity $\hat{g}(\bm{Q} \bm{F}) = \hat{g}(\bm{F})$ with respect to
$\bm{F}$ yields the transformation rule of the gradient,
\begin{equation}
D\hat{g} (\bm{Q}\bm{F})
 = \bm{Q}\,D\hat{g} (\bm{F})\,,
\label{eq-gradgObj}
\end{equation}
and therefore each of the three terms in Eq.~\eqref{eq-bregman} is invariant:
the first two trivially, and the third because
\begin{equation}
D\hat{g} (\bm{Q}\bm{F}) : (\bm{Q}\,\nabla_0 \bm{\varphi} - \bm{Q}\bm{F})
 = \bm{Q} D\hat{g} (\bm{F}) : \bm{Q}\,(\nabla_0 \bm{\varphi} - \bm{F})
 = D\hat{g} (\bm{F}) : (\nabla_0 \bm{\varphi} - \bm{F})\,,
\end{equation}
where the last equality follows from $\bm{Q}\bm{A} : \bm{Q}\bm{B} = \bm{A} : \bm{Q}^{\mathsf T}\bm{Q}\,\bm{B} = \bm{A} : \bm{B}$. Hence $D_{B}(\hat{g},\bm{Q}\,\nabla_0 \bm{\varphi} , \bm{Q}\bm{F}) = D_{B}(\hat{g},\nabla_0 \bm{\varphi},\bm{F})$ for every rotation
$\bm{Q}$, provided that $\hat{g}$ is itself objective. 

Employing the proposed definition of the stabilizing function, we can calculate the
variational derivatives of the functional~\eqref{eq-HWenergy-large} leading to the mixed
variational problem: find $(\bm{\varphi},\bm{F},\bm{P}) \in \mathcal{X}\times\mathcal{F}\times\mathcal{F}$ such that
\begin{equation}
    \begin{split}
        \int_{\mathcal{B}_0} {\bm{P}:\nabla_0 \delta\bm{\varphi}} \ dV + \alpha \int_{\mathcal{B}_0} { \left[ D\hat{g}(\nabla_0 \bm{\varphi}) - D\hat{g}(\bm{F}) \right] :\nabla_0 \delta\bm{\varphi}} \ dV
        - \int_{\mathcal{B}_0} {\bm{b}_0 \cdot \delta\bm{\varphi}} \ dV - \int_{{\partial_t\mathcal{B}}_0}{\bm{t}_0 \cdot \delta\bm{\varphi}} \ dA &= 0 \quad , \\
        \int_{\mathcal{B}_0} {\left(D\phi({\bm{F}}) - \bm{P}\right): \delta\bm{F}} \ dV - \alpha \int_{\mathcal{B}_0} {\left[D^2\hat{g}(\bm{F}): (\nabla_0 \bm{\varphi}-\bm{F}) \right]:\delta\bm{F}} \ dV &= 0 \quad , \\
        \int_{\mathcal{B}_0} {(\nabla_0 \bm{\varphi}-\bm{F}):\delta\bm{P}} \ dV &= 0 \quad , 
    \end{split}
    \label{eq-huwashizu-LS}
\end{equation}
for all $(\delta\bm{\varphi},\delta\bm{F},\delta\bm{P}) \in \mathcal{X}_0\times\mathcal{F}\times\mathcal{F}$. Notably, the choice~\eqref{eq-bregman} ensures that requirements R2 and R3 are automatically satisfied.

\subsection{Finite element discretization}
\label{sect-mixedformprob-largeFE}

To introduce the discrete formulation of the problem, we assume that the reference configuration can
be partitioned into primal and dual meshes, just like the domain of the small strain problem in
Section~\ref{sect-review-averaging}. Then, we define the spaces
\begin{equation}
    \begin{split}
        \mathcal{X}^h &= \{ \bm{\varphi}^h \in \mathcal{X} \ , \bm{\varphi}^h|_{e_{0\alpha}} \in [R_1(e_{0\alpha})]^d \ , e_{0\alpha} \in \mathcal{M} \} \ , \\
        \mathcal{X}_0^h &=\{ \bm{\varphi}^h \in \mathcal{X}_0 \ , \bm{\varphi}^h|_{e_{0\alpha}} \in [R_1(e_{0\alpha})]^d \ , e_{0\alpha} \in \mathcal{M} \} \ , \\
        \mathcal{F}^h &=\{ \bm{F}^h \in \mathcal{F} \ , \bm{F}^h|_{{e_0}^\beta} \in [R_0({e_0}^\beta)]^{d\times d} \ , {e_0}^\beta \in {\mathcal{M}^*} \} \ ,
    \end{split}
    \label{eq-spacesFE-LS}
\end{equation}
where $e_{0\alpha}$ are the elements of the primal mesh and ${e_0}^\beta$ are the cells of the dual mesh, both in reference configuration. 

The finite element formulation is obtained as the discrete counterpart of
problem~\eqref{eq-huwashizu-LS}. Following the same steps of Section~\ref{sect-review-averaging},
a deformation finite element problem can be obtained: find $\bm{\varphi}^h \in \mathcal{X}^h$ such that
\begin{equation}
\begin{split}
    \int_{\mathcal{B}_0} {\bm{P}^h(\bm{\varphi}^h):\nablabar_0 \delta\bm{\varphi}^h} \ dV - \alpha \int_{\mathcal{B}_0} {D\hat{g}\left(\nablabar_0 \bm{\varphi}^h\right):\nablabar_0 \delta\bm{\varphi}^h} \ dV + \alpha \int_{\mathcal{B}_0} {D\hat{g}(\nabla_0 \bm{\varphi}^h):\nabla_0 \delta\bm{\varphi}^h} \ dV = \\
    = \int_{\mathcal{B}_0} {\bm{b}_0 \cdot \delta\bm{\varphi}^h} \ dV + \int_{{\partial_t\mathcal{B}}_0}{\bm{t}_0 \cdot \delta\bm{\varphi}^h} \ dA  \quad , \quad \forall \delta\bm{\varphi}^h \in \mathcal{X}_0^h \ ,
\end{split}
\label{eq-nodalFEM-LS}
\end{equation}
where
\begin{equation}
    \bm{F}^h(\bm{\varphi}^h) := \nablabar_0 \bm{\varphi}^h := \sum_{\beta=1}^{N_n} \bm{F}^h_\beta \ \chi^\beta \quad , \quad \bm{F}^h_\beta := \left. \nablabar_0 \bm{\varphi}^h \right|_{e_0^\beta} := \frac{1}{|{e_0}^\beta|} \int_{{e_0}^\beta} {\nabla_0 \bm{\varphi}^h} \ dV \ ,
    \label{eq-defgradfield}
\end{equation}
\begin{equation}
    \bm{P}^h(\bm{\varphi}^h) := \sum_{\beta=1}^{N_n} \bm{P}^h_\beta \ \chi^\beta \quad , \quad \bm{P}^h_\beta := \frac{1}{|{e_0}^\beta|} \int_{{e_0}^\beta} {D\phi({\bm{F}^h})} \ dV = D\phi({\bm{F}^h_\beta}) \ .
    \label{eq-piolaKfield}
\end{equation}
Note the expected separation of the stabilization term in purely dual-mesh and primal-mesh integrals
in Eq.~\eqref{eq-nodalFEM-LS} and the unaltered constitutive relationship in
Eq.~\eqref{eq-piolaKfield}.

As for small strain kinematics, the assumed deformation gradient~\eqref{eq-defgradfield} can be interpreted as the result
of a modified gradient operator $\nablabar_0$ on the discrete deformation.
%
%

\begin{remark}
   As is well known, balance equations in the large-strain setting can be expressed in terms of different work-conjugate pairs. In the case of the nodal strain method, stating the Hu-Washizu functional in terms of different strain and stress measures in the continuum problem results in different discrete approximations, which means that one must choose the main kinematic descriptor for the nodal approach. This happens because the projection operator applied to one strain measure does not carry over to other measures due to the nonlinear relations between them. As an example, consider the computation of the Green-Lagrange strain tensor $\bm{E}(\bm{G}) = 1/2 \ (\bm{G}^T \bm{G} - \bm{1})$ for a certain deformation gradient $\bm{G}$. Then, given the description~\eqref{eq-defgradfield}, it is clear that in general
\begin{equation}
    \bm{E}(\nablabar_0 \bm{\varphi}^h) \neq {\overline{\bm{E}}^h(\nabla_0 \bm{\varphi})} \ .
\end{equation}
\end{remark}

\begin{remark}
The structure of the stabilization term deserves a detailed
discussion, since it constitutes one of the contributions of this work. Requirement R2 stated at the beginning of this section is often associated, by analogy with the small strain case stabilization term~\eqref{eq-stabss}, with a function $\hat{\phi}({\cdot})$ quadratic in the difference of the kinematic descriptors, see for instance \cite{puso2006stabilized}. Interestingly, the Bregman form $D_{B}$ of Eq.~\eqref{eq-bregman} generalizes this idea for an arbitrary smooth $\hat{g}$, because (a)~its $\nabla_0 \bm{\varphi}$-variation splits into the two single-argument terms of the first equation in Eq.~\eqref{eq-huwashizu-LS}; and (b)~its $\bm{F}$-variation is proportional to the intra-cell fluctuation $\nabla_0 \bm{\varphi}^h - \bm{F}^h$, which vanishes due to the projection operation. In fact, on the discrete spaces the stabilization energy takes the remarkably simple form

\begin{equation}
\alpha \int_{\mathcal{B}_0} {D_{B} (\hat{g},\nabla_0 \bm{\varphi}^h,\bm{F}^h)} \ dV
 = \alpha \int_{\mathcal{B}_0} {\hat{g}(\nabla_0 \bm{\varphi}^h)} \ dV 
 - \alpha \int_{\mathcal{B}_0} {\hat{g}(\bm{F}^h)} \ dV \,,
\label{eq-energydiscrete}
\end{equation}
since the linear correction term in Eq.~\eqref{eq-bregman} integrates to zero
cell by cell in the dual mesh.   
\end{remark}

\begin{remark}
Expression~\eqref{eq-energydiscrete} allows a clear interpretation of the function $\hat{g}$: this function \emph{resembles} an energy density, and the stabilization consists of the difference between the stabilizing energy sampled on the primal mesh and the same energy sampled on the dual mesh.

This is precisely the structure of the energy-sampling stabilization proposed by \cite{sivapuram2019energy} for small strains, which was introduced there as a numerical device. The present derivation shows that it descends from the augmented Hu-Washizu functional~\eqref{eq-HWenergy-large}, provided the stabilizing term is understood as the Bregman divergence of the sampled energy. In the same spirit, this result clarifies the status of stabilization strategies that reuse the constitutive model of the real material, or part of it, as stabilizing material \cite{gee2009uniform}: such a choice can be made consistent with an energy functional like~\eqref{eq-HWenergy-large}, but only through the Bregman construction~\eqref{eq-bregman}, and not by direct evaluation of the material energy on the difference of the kinematic descriptors.

The admissible class of stabilizers compatible with requirements R1, R2, and R3 is therefore identified
as the family of Bregman divergences of frame-indifferent functions $\hat{g}$. Further discussion around an effective choice of $\hat{g}$ is presented in Section~\ref{sect-stablarge}.
\end{remark}

\subsection{The projection operator in a change of reference configuration}
\label{sect-projop-reference}

In a finite element simulation of a finite strain problem, each time the body is remeshed the new
mesh becomes the reference configuration for the Lagrangian elements and the definition of material
derivatives must be updated accordingly. As a result, the assumed gradient
operator~\eqref{eq-defgradfield} that takes care of the averaging operation must also be updated.

Consider, in addition to the initial configuration ${\mathcal{B}_0}$ and the current configuration $\mathcal{B}$, an arbitrary
intermediate configuration $\mathcal{B}_I$. In this case, two new deformation maps can be introduced: $\bm{\varphi}_0 : {\mathcal{B}_0} \rightarrow
\mathcal{B}_I$ and $\bm{\varphi}_I : \mathcal{B}_I \rightarrow \mathcal{B}$ such that
$\bm{\varphi}=\bm{\varphi}_I\circ\bm{\varphi}_0$. Thus, the deformation gradient can be written as
\begin{equation}
  \label{eq-change-ref-cont}
    \nabla_0 \bm{\varphi} = \nabla_{I} \bm{\varphi}_I \cdot \nabla_0 \bm{\varphi}_0 \ ,
\end{equation}
where $\nabla_I$ is the gradient with respect to the intermediate configuration.
When employing a finite element discretization, we define $\bm{\varphi}^h_0$ and
$\bm{\varphi}^h_I$ to be, respectively, the discrete deformation at the moment of the last
remeshing and the deformation from that configuration to the current one. In other words, the intermediate configuration $\mathcal{B}_I$ always coincides with the configuration of the last remeshing. Then,
Eq.~\eqref{eq-change-ref-cont} takes the discrete representation
\begin{equation}
    \nabla_0 \bm{\varphi}^h = \nabla_{I} \bm{\varphi}^h_I  \cdot\nabla_0 \bm{\varphi}_0^h \ .
    \label{eq-refconfchange}
\end{equation}

Let $\mathcal{B}_I$ be the last chosen reference configuration, then $\nabla_0 \bm{\varphi}_0^h$ is a
known and fixed field for the subsequent computations. This object captures the (local) deformation
from the initial configuration to the one where the mesh is defined for the last time, and it should remain clear that this is a quantity computed in the primal mesh. Then, replacing expression~\eqref{eq-refconfchange} in the nodal representation~\eqref{eq-defgradfield}
\begin{equation}
    \bm{F}^h_\beta := \left. \nablabar_0 \bm{\varphi}^h \right|_{e_0^\beta} := \frac{1}{|{e_0}^\beta|} \int_{{e_0}^\beta} {\nabla_I \bm{\varphi}^h_I  \cdot \nabla_0 \bm{\varphi}_0^h} \ dV  \quad , \quad |{e_0}^\beta| := \int_{{e_0}^\beta} {} \ dV
    \label{eq-confchange-Fbar1}
\end{equation}


Finally, changing the configuration using well-known rules in continuum mechanics, and using the notation $J_0 := \det(\nabla_0 \bm{\varphi}_0^h)$ it is possible to write the bar operator in terms of the dual cells ${e_I}^\beta$ in the last reference configuration $\mathcal{B}_I$ as
\begin{equation}
    \left. \nablabar_0 \bm{\varphi}^h \right|_{e_0^\beta} := \frac{1}{|{e_0}^\beta|} \int_{{e_I}^\beta}{\nabla_I \bm{\varphi}^h_I  \cdot \nabla_0 \bm{\varphi}_0^h \ J_0^{-1} \ dV}  \quad , \quad     |{e_0}^\beta| = \int_{{e_I}^\beta}{J_0^{-1} \ dV} \ .
\label{eq-confchange-Fbar2}
\end{equation}



Expression~\eqref{eq-confchange-Fbar2} requires the tracking of $\nabla_0 \bm{\varphi}_0^h$, which must be computed every time a change in reference configuration happens. This quantity is easily computed using the standard finite element functions of the primal mesh. Note that the bar operator~\eqref{eq-confchange-Fbar2} takes gradients with respect to the last intermediate configuration $\mathcal{B}_I$ and returns the nodal gradients with respect to the initial configuration $\mathcal{B}_0$. Therefore, the weak form~\eqref{eq-nodalFEM-LS} is still computed in the initial configuration.

\subsection{On the choice of the stabilizing function $\hat{g}$}
\label{sect-stablarge}

Section~\ref{sect-review-stabilization} stressed the need to include a
stabilization term in the Hu-Washizu functional for mixed finite elements in
the small strain problem. Naturally, the
difficulties increase in the large-strain case, and, although there are no theoretical proofs
showing the well-posedness of the weak form in the nonlinear case, this section illustrates the issues related to an incorrectly chosen stabilization term~\eqref{eq-bregman}, 
and presents a discussion of the role of the function $\hat{g}$.

First, note that it is not necessary for $\hat{g}$ to be a stored energy function. 
This is so because when the difference $\nabla_0 \bm{\varphi}^h-\bm{F}^h\to 0$ due to
mesh refinement, the stabilization~\eqref{eq-energydiscrete} added in the discrete functional
vanishes, independently of any physical interpretation. However, from a practical point of view and
dimensional considerations, it is convenient that some \emph{energy-like} features remain present
in the stabilization.

To motivate this remark, note that inspection of the average operation~\eqref{eq-defgradfield} reveals that a well-behaved deformation gradient $\bm{F}^h$ in the dual mesh, that is, one with $\det(\bm{F}^h)>0$, does not imply the same condition for its primal mesh counterpart $\nabla_0 \bm{\varphi}^h$. The domain of $\hat{g}$ is part of the modelling choice and controls the behaviour of the stabilization under severe element distortion. If $\hat{g}$ remains finite as $\det(\nabla_0 \bm{\varphi}^h) \to 0$ on an element of the primal mesh, then the stabilization energy of a configuration containing degenerate, or even inverted, elements is finite, and nothing prevents such states from being reached during the solution process. Such a situation is especially problematic for guaranteeing the reliability of a remeshing procedure based on flips of groups of elements of the discarded mesh, as is the case in this paper. To avoid this type of behaviour, the choice of $\hat{g}$ must satisfy the growth condition $\hat{g}(\bm{G}) \to +\infty$ as $\det\bm{G} \to 0^+$, characteristic of realistic strain energy densities. With this choice, the primal-mesh contribution in Eq.~\eqref{eq-energydiscrete} acts as a barrier against element degeneration and inversion. The influence of $\hat{g}$ in the solution is exemplified below.

\subsubsection{Numerical example}
\label{subs-stability-large-example}
\begin{figure}[t]
    \centering
    \includegraphics[width=0.26\textwidth]{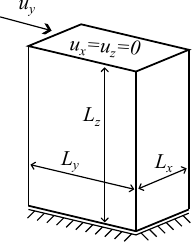}
    \caption{Problem setting.}
    \label{fig-stabEx-LS}
\end{figure}

As in Section~\ref{subs-stability-small-example}, we use a numerical example to illustrate the need
for a correct stabilization term in the finite strain, nodal finite element method. We consider a solid
brick of dimensions $L_x$, $L_y = 2L_x$, $L_z=3L_x$, with the whole displacement field constrained on its
bottom face, and  imposed displacements on its top face in the $y$ direction of $u_y = L_y$; see
Figure~\ref{fig-stabEx-LS}. The body is modelled as elastic with a Saint Venant-Kirchhoff constitutive law, with Young's modulus $E = 1$ and Poisson's ratio $\nu = 0.4$, or equivalently with Lam\'e parameters $(\lambda,\mu) = (2,1/2) \, 5/7$.

Three cases are compared: a nodally integrated finite element method without stabilization, followed by two stabilized cases with stabilization functions $\hat{g}_1$ and 
$\hat{g}_2$, respectively. Consider first the quadratic function
\begin{equation}
\hat{g}_1(\bm{G}) := \frac{\hat{\mu}}{2}\|\bm{G}\|^2\,, \qquad \|\bm{G}\|^2 := \text{tr}(\bm{G}^\text{T}\bm{G}) \,,
\label{eq-ghat1}
\end{equation}
This choice is included to emphasize that $\hat{g}_1$ is clearly not a realistic energy density, since neither the function nor its derivative with respect to $\bm{G}^\text{T}\bm{G}$ (which would result in the associated second Piola-Kirchhoff stress) vanish for $\bm{G} = \bm{1}$. However, it does verify the objectivity requirement stated in Section~\ref{sect-mixedformprob-large-stab}. Note that for the choice~\eqref{eq-ghat1} the stabilizing function is 
\begin{equation}
    \hat{\phi}({\nabla_0 \bm{\varphi},\bm{F}}) = D_{B}(\hat{g},\nabla_0 \bm{\varphi},\bm{F}) = \frac{\hat{\mu}}{2} ||\nabla_0 \bm{\varphi}-\bm{F}||^2 \, ,
\end{equation}
close to a proposal from literature \cite{puso2006stabilized}. Consider as a second choice the compressible neo-Hookean density
\begin{equation}
\hat{g}_2(\bm{G}) := \frac{\hat{\mu}}{2}\left(\|\bm{G}\|^2 - 3\right)
             - \hat{\mu}\,\ln J
             + \frac{\hat{\lambda}}{2}\,(\ln J)^2\,,
\qquad J := \det\bm{G} \,,
\label{eq-ghat2}
\end{equation}
whose volumetric barrier penalizes the degeneration of primal-mesh elements. Note that the real material model, a Saint Venant-Kirchhoff model, does not correspond to any of the two stabilizing functions, which are defined independently, and the problem remains variationally consistent. The specific values adopted in each case are summarized in Table~\ref{table-stabexample}.

\begin{table}[ht]
\centering
    \begin{tabular}{c c c}
         parameters & $\hat{g}_1$, Eq.~\eqref{eq-ghat1} & $\hat{g}_2$, Eq.~\eqref{eq-ghat2} \\
        \hline
         $\alpha$ & 0.5 & 0.1 \\
         $\hat{\mu}$ & $\mu$ & $\mu$ \\
         $\hat{\lambda}$ & - & 0.0
    \end{tabular}
    \caption{Stabilization coefficient and material parameters for stabilized formulations}
    \label{table-stabexample}
\end{table}

\begin{figure}[htbp]
    \centering

    \begin{subfigure}{\textwidth}
        \centering
        \includegraphics[width=0.9\linewidth]{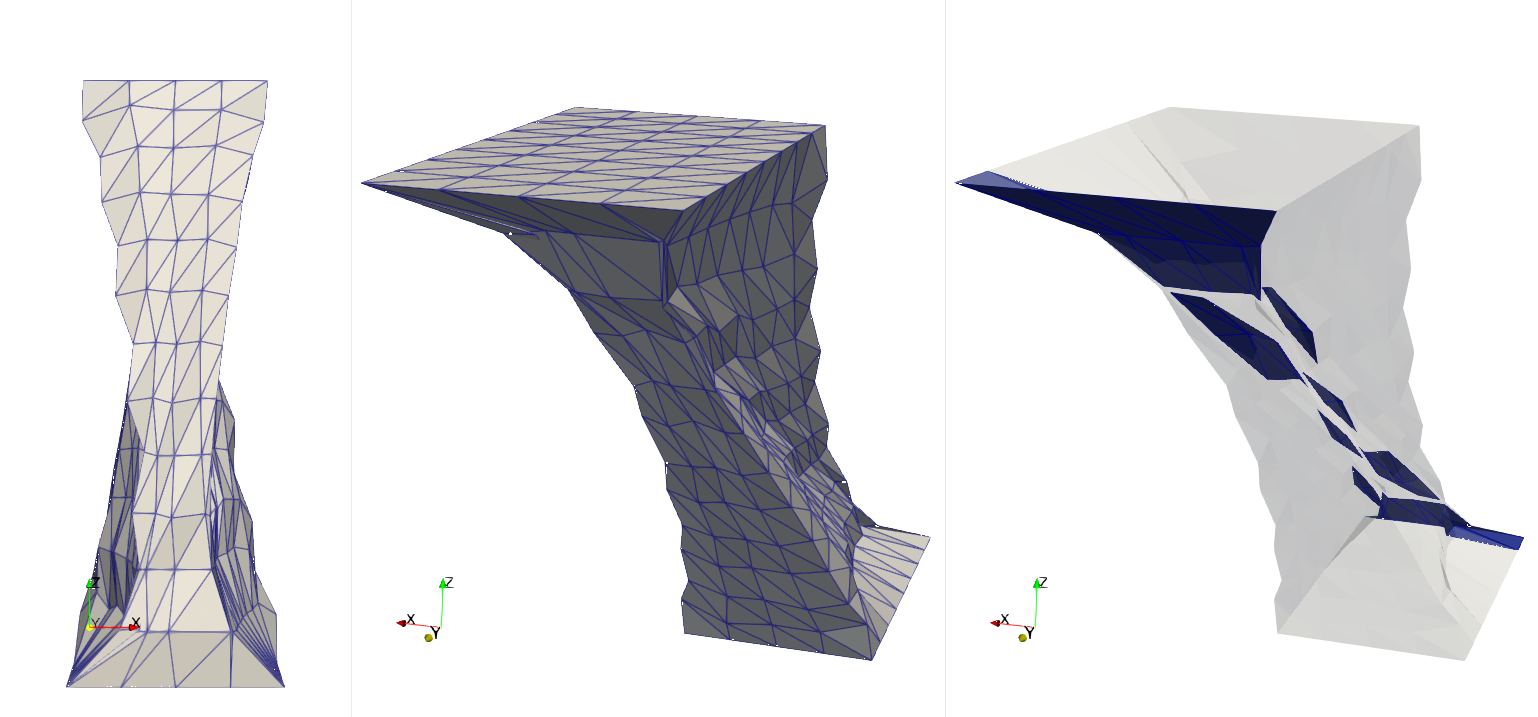}
        \caption{Deformed configuration with no stabilization.}
        \label{fig-stablarge-without}
    \end{subfigure}
    \vspace{0.1cm}

    \begin{subfigure}{\textwidth}
        \centering
        \includegraphics[width=0.9\linewidth]{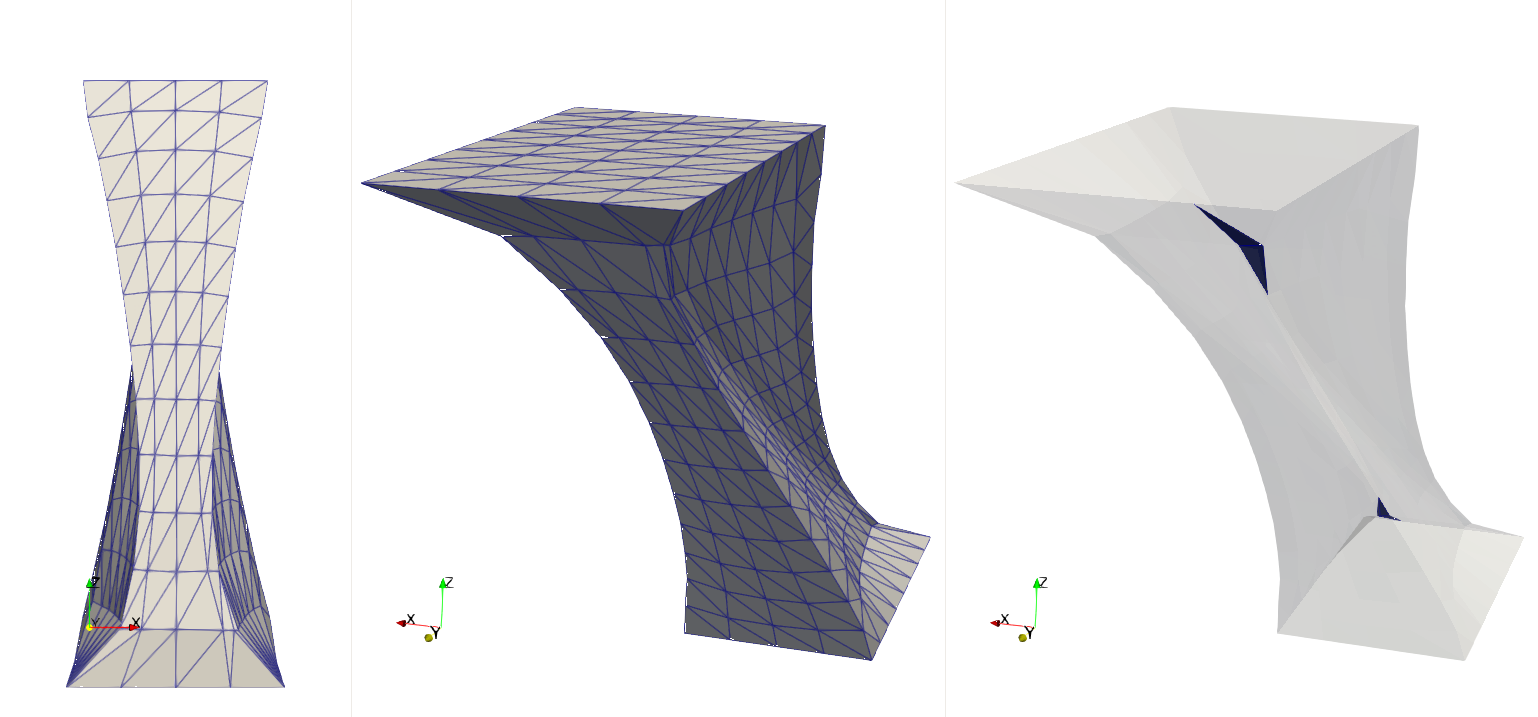}
        \caption{Deformed configuration with stabilizing function~\eqref{eq-ghat1}.}
        \label{fig-stablarge-with-fake}
    \end{subfigure}
    \vspace{0.1cm}

    \begin{subfigure}{\textwidth}
        \centering
        \includegraphics[width=0.9\linewidth]{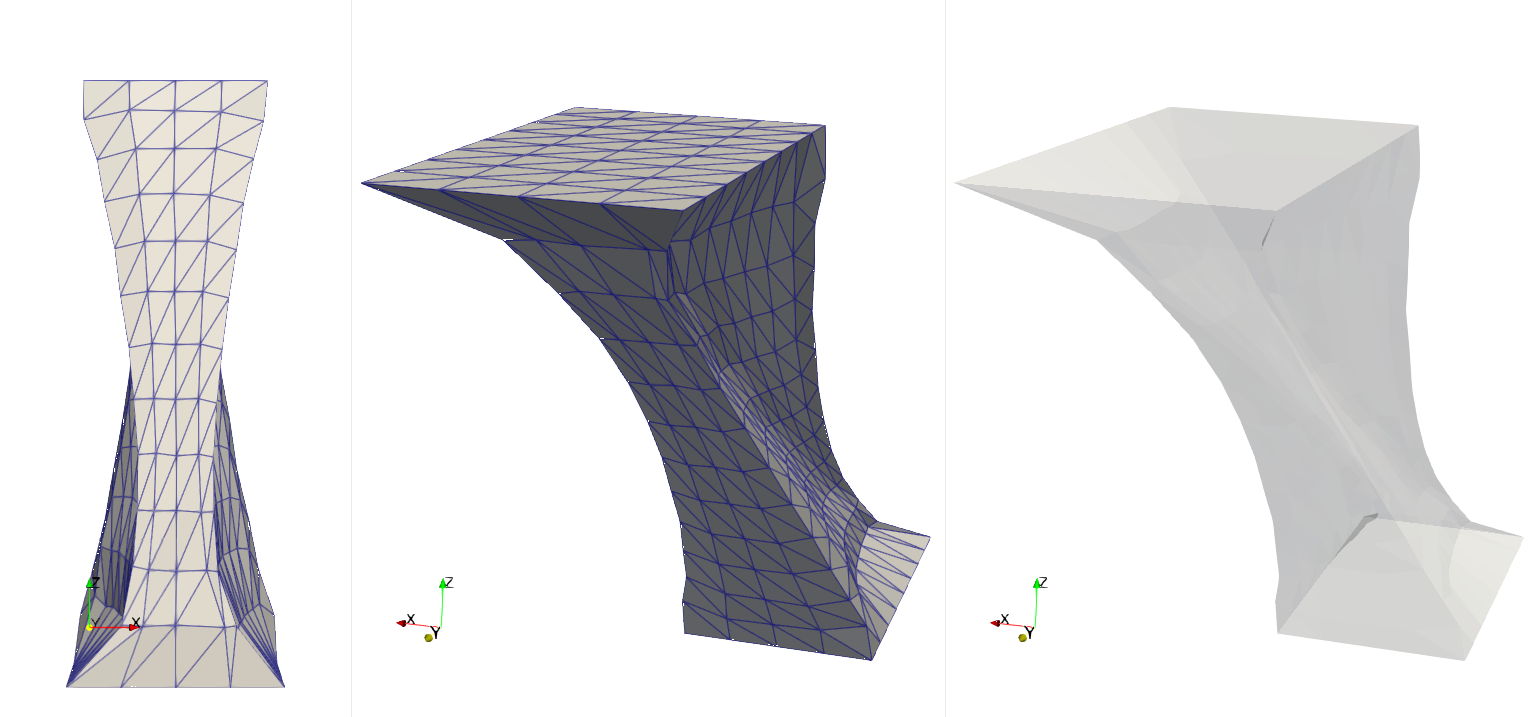}
        \caption{Deformed configuration with stabilizing function~\eqref{eq-ghat2}.}
        \label{fig-stablarge-with-nh}
    \end{subfigure}

    \caption{Comparison of deformed configurations of the nodal strain approach in large-strains regime without and with stabilization term included. Left and centre: different views of the deformed meshes. Right: inverted elements in dark blue inside the transparent body.}
    \label{fig-stabexample-large}
\end{figure}

Figure~\ref{fig-stabexample-large} compares the deformed configurations obtained for the three
cases, where the left and centre figures show different views of the deformed primal mesh, and the
right figures show in dark colour the inverted elements of the primal mesh. In
Figure~\ref{fig-stabexample-large}(a), an oscillatory deformation pattern can be observed for the non-stabilized formulation. In this solution, elements of both the primal and the dual meshes attained
negative values of the determinant of their respective deformation gradient fields.
Figure~\ref{fig-stabexample-large}(b) shows how the stabilizing function~\eqref{eq-ghat1}
successfully removes the oscillations from the solution. However, due to the lack of a barrier
against element degeneration, some elements of the primal mesh still get inverted. Finally,
Figure~\ref{fig-stabexample-large}(c) illustrates how a properly defined stabilization function
achieves both goals of stabilizing and avoiding element inversion at the same time, despite the
smaller stabilizing coefficient in comparison with case~(b).

\section{Application to history-dependent problems with remeshing}
\label{sect-remesh-diff}

This section explains how the nodal structure of the proposed formulation can be exploited to avoid
remapping of internal variables when remeshing. As a
result, this property will be shown to reduce the diffusion of information in the history variables
that model the inelastic processes. In this sense, the section collects the three ingredients
required by the numerical study performed in Section~\ref{sect-examples}: Section~\ref{sect-epmodel}
presents the history-dependent constitutive model whose internal variables must be preserved,
Section~\ref{sect-remeshing} describes the algorithmic procedure and what a typical remeshing step
implies for the nodal approximation, and Section~\ref{sect-diffusion} introduces an indicator used to quantify the level of diffusion of a given field.

\subsection{Finite-strain elastoplastic constitutive model}
\label{sect-epmodel}

The history-dependent response considered in this work is described by the classical finite strain
$J_2$ plasticity model formulated in principal logarithmic
stretches~\cite{simo1992algorithms,simo1998numerical}. Since the model is standard, only its
essential ingredients are summarized here; the implementation employed is that of the open-source material library MUESLI~\cite{portillo2017muesli}.

The model assumes that the material state is fully described by the set
$\{\bm{F},\bm{b}^e, \xi, \boldsymbol{\xi}\}$. In this set, $\bm{F}$ is the deformation gradient and admits the multiplicative decomposition
$\bm{F} = \bm{F}^{e}\bm{F}^{p}$, into elastic and plastic parts, with plastically incompressible flow,
$\det\bm{F}^{p} = 1$. The other three terms define the set of internal variables. Specifically,
$\bm{b}^e := \bm{F}^{e}\bm{F}^{e\text{T}} = \bm{F}\, (\bm{F}^{p\text{T}}\bm{F}^{p})^{-1} \bm{F}^{\text{T}}$ is the elastic left Cauchy-Green tensor,
$\xi$ is the isotropic hardening strain-like variable, and $\boldsymbol{\xi}$ is the kinematic hardening strain-like
variable. Note that $\xi$ is also identified as an equivalent plastic strain.

The thermodynamics of a material point is governed by a free energy of the form
\begin{equation}
    \phi = \phi({\bm{b}^e, \xi, \boldsymbol{\xi}}) := W(\bm{b}^e) + \mathcal{H}(\xi, \boldsymbol{\xi}) \;,
\end{equation}
where it is assumed that the free energy can be additively decomposed into an elastic stored energy function~$W$ and a potential function~$\mathcal{H}$ for the hardening variables.

The elastic response is hyperelastic of Hencky type. Denote by $(\lambda^{e}_{i})^2$ and $\bm{n}_{i}$ the eigenvalues and eigenvectors of
$\bm{b}^e$, and by $\varepsilon^{e}_{i} := \tfrac{1}{2}\ln\lambda^{e}_{i}$ the
principal logarithmic elastic strains. The Hencky stored energy is quadratic in the
logarithmic elastic strains $\boldsymbol{\varepsilon}^{e} := \sum_{i} \varepsilon^{e}_{i}\, \bm{n}_{i}\otimes\bm{n}_{i}$, 
\begin{equation}
W(\bm{b}^e(\boldsymbol{\varepsilon}^{e}))
 := \frac{\lambda}{2}\big(\text{tr}\,\boldsymbol{\varepsilon}^{e}\big)^{2}
 + \mu\,\|\boldsymbol{\varepsilon}^{e}\|^{2} \,,
\label{eq-hencky-energy}
\end{equation}
where $(\lambda,\mu)$ are the Lam\'e constants. In turn, the hardening potential is quadratic in the strain-like internal variables
\begin{equation}
\mathcal{H}(\xi, \boldsymbol{\xi}) := \frac{1}{2} K \xi^2 + \frac{1}{3}H \|\boldsymbol{\xi}\|^{2}\,,
\end{equation}
where $K$ and $H$ are the isotropic and kinematic hardening moduli. From standard arguments, it is possible to derive the (Kirchhoff) stress tensor $\boldsymbol{\tau}$, and conjugate forces~$q$ and~$\bm q$ from the free energy
\begin{equation}
    \boldsymbol{\tau} := \partial_{\boldsymbol{\varepsilon}^{e}}\phi = \partial_{\boldsymbol{\varepsilon}^{e}} W \; , \qquad q := -\partial_\xi \phi = -\partial_\xi \mathcal{H} \; , \qquad \bm q := -\partial_{\boldsymbol{\xi}} \phi = -\partial_{\boldsymbol{\xi}} \mathcal{H} \; .
\end{equation}

The elastic domain is defined by the von Mises yield criterion in principal Kirchhoff stress space,
\begin{equation}
f(\boldsymbol{\tau}, q, \bm q)
 := \big\| \text{dev}[\boldsymbol{\tau}] - \bm q \big\|
 - \sqrt{\frac{2}{3}}\,\left(\sigma_{Y} - q\right) \;\leq\; 0\,,
\label{eq-yield}
\end{equation}
where $\sigma_{Y}$ is the yield stress. The evolution equations are associative, obtained from the principle of maximum dissipation. The numerical integration of the model follows the exponential return mapping
algorithm explained in \cite{simo1992algorithms, simo1998numerical}.

Within the nodal finite element formulation, the elastoplastic model is evaluated once per dual cell: the internal variables belong to the same piecewise-constant space over the dual mesh as the deformation gradient, and the update is driven by the nodal averaged gradients $\bm{F}^h_\beta$, so that the nodal stress in the discrete problem is the output of the algorithmic constitutive function, generalizing the hyperelastic relation of the previous section to the inelastic setting. 

\subsection{Remeshing}
\label{sect-remeshing}

As previously discussed, the nodal strain finite element method developed in this work behaves as a
particle method in the sense that the complete material information can be associated with the nodes,
which behave not only as carriers of the kinematics, but also as material points. Therefore, when
element connectivity is redefined during remeshing, remapping of the internal variables is unnecessary
since node positions remain unchanged. This is in contrast with standard FEM formulations where
changes in the location of the Gauss points require the interpolation or projection of the internal
variables, leading to diffusion of information. To isolate and quantify
the performance of the method to reduce diffusion due to transfer of information, all the examples
presented in this work limit remeshing to the reconnection of nodes, or in other words redefinition
of elements by flip operations~\cite{freitag1997tetrahedral,rivara1997ws}. 

\begin{figure}[t]
    \centering
    \includegraphics[width=\textwidth]{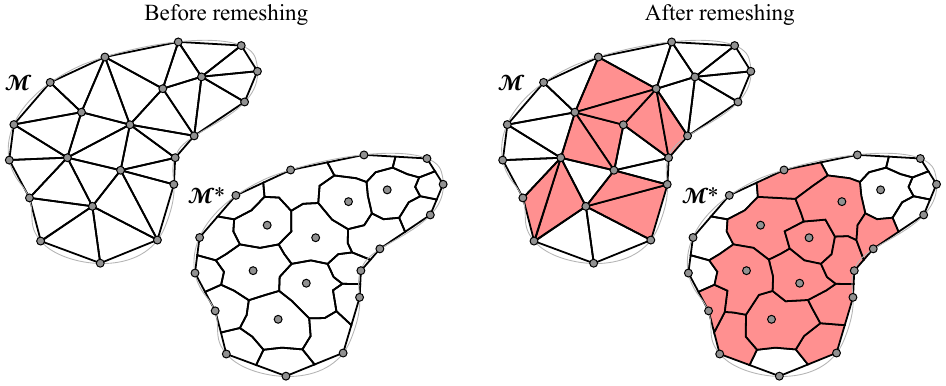}
    \caption{Change in dual mesh due to remeshing of the primal mesh. Highlighted elements changed due to the remeshing.}
    \label{fig-remesh}
\end{figure}

Despite the absence of history remapping, redefinition of the element connectivity in the primal
mesh requires the redefinition of the dual mesh. This in turn means that the function spaces of the
discrete problem~\eqref{eq-spacesFE-LS} are changed during remeshing. These modifications perturb
the equilibrium and, therefore, the solution must be rebalanced after the redefinition of the dual
mesh. Figure~\ref{fig-remesh} illustrates this process. As shown, modifications in the connectivity in the
primal mesh induce changes in the dual mesh.

Furthermore, due to the large distortions that a deforming body may suffer, in general, it is necessary to perform a reference configuration update together with the remeshing step. This is because the new elements, created in the current configuration, and their shape functions may not be correctly defined in the initial configuration. In this scenario the bar operator to compute the gradients of the shape functions must be modified following the details explained in Section~\ref{sect-projop-reference}. Note however that the deformation gradient $\nabla_0 \bm{\varphi}_0^h$ is a quantity of the primal mesh, and needs, therefore, to be remapped between old and new primal meshes. In other words, some diffusion of information exists in the nodal strain approach, although it is not related to the material history variables. The diffusion of $\nabla_0 \bm{\varphi}_0^h$, also present in standard finite elements, can be reduced following the ideas discussed in \cite{rodriguez2016particle, rodriguez2017continuous}.

\subsection{Measuring diffusion}
\label{sect-diffusion}

Standard remeshing strategies of finite element inelastic models \emph{smear} the
information carried by the internal variables at the material points. There is no canonical way to
account for this ``information diffusion'' that spreads fields such as the plastic slip,
damage, etc. Since we would like to assess the effects of remeshing on discretization schemes, we
propose next a measure that will be used in the numerical examples of Section~\ref{sect-examples}.

First, we note that the capability of a field to diffuse in a specific domain is related to how
close this field is to a homogeneous distribution. However, quantifying the degree of homogeneity
or, on the contrary, locality of a certain field is not trivial. For instance, one could think
of the measure of volume where this field is greater than zero as a possible metric. This, however,
has both fundamental and implementation-related drawbacks. First, even if the field is positive in the whole domain, it is capable of diffusion if it is not homogeneous. Furthermore, this measure would be very sensitive to small numerical errors due to floating point accuracy.

In this work, a measure is adopted based on concepts borrowed from information theory. Let $p$ be a probability density function defined on the domain $\mathcal{B}$, with the property $\int_\mathcal{B} {p} \ dV=1$. The set $\mathcal{B}_p$ where $p>0$ is referred to as the support of $p$; clearly, $\mathcal{B}_p \subseteq \mathcal{B}$. Now, let $s$ be a certain non-negative scalar field defined on the domain, then it is possible to define the probability distribution associated with the field~$s$ as
\begin{equation}
    p(\bm{X}) := \frac{s(\bm{X})}{\int_\mathcal{B} {s} \ dV} \ ,
    \label{eq-probdistr}
\end{equation}
assuming the denominator is positive. The case of ${\int_\mathcal{B} {s} \ dV = 0}$ is addressed later in this section. Then, given a subset $\mathcal{A}\subseteq \mathcal{B}$, the value $\int_\mathcal{A} p \ dV$ is the contribution from the region $\mathcal{A}$ to the total $s$-content in the whole domain. Additionally, let $q$ be a \emph{reference} probability distribution also defined on $\mathcal{B}$. It is possible to get an idea of how close $p$ is to $q$ by means of the \emph{relative entropy} or Kullback–Leibler divergence $D(p||q)$, defined as 
\begin{equation}
    D(p||q) := \int_\mathcal{B} {p \ln{\frac{p}{q}}} \ dV \ ,
    \label{eq-KLdivergence}
\end{equation}
with the convention $0 \ln{0} = 0$. The relative entropy has the useful property of being always nonnegative, and zero if and only if $p = q$, almost everywhere. Naturally, the reference distribution used to measure proximity to the homogeneous case is:
\begin{equation}
    q_h := \frac{1}{|\mathcal{B}|} \ .
\end{equation}

Then, the following indicator is adopted as a measure of how diffused (or how close to homogeneous) a field is

\begin{equation}
    f_{eq} := \exp\{-D(p||q_h)\} \ .
    \label{eq-diffindicator}
\end{equation}

Some properties of this indicator are:

\begin{itemize}
    \item It takes values only in $(0,1]$. Values close to $0$ indicate that the field $s$
      is very localized, while $f_{eq}=1$ represents the homogeneous case.
      
    \item Due to the exponential law, the indicator is more
      sensitive to changes in the relative entropy when $p$ is close to the homogeneous case $q_h$ than
      when it is highly concentrated. This, in turn, makes it insensitive to machine-precision zeros
      that may appear in initial stages.
      
    \item Consider the probability distribution defined by $\hat{p}:=1/|\Hat{\mathcal{B}}|$ if $x\in\Hat{\mathcal{B}}$. That is, a distribution that is homogeneous in a subdomain $\Hat{\mathcal{B}} \subset \mathcal{B}$. It can be proven that in this situation $f_{eq} = |\Hat{\mathcal{B}}|/|\mathcal{B}|$.
\end{itemize}

\begin{figure}[t]
    \centering
    \includegraphics[width=12cm]{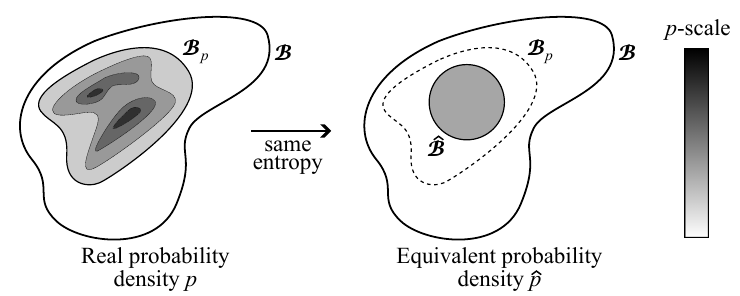}
    \caption{Scheme of the entropy equivalence between probability densities $p$ and $\hat{p}$, that allows the interpretation of $f_{eq}$ as a volume fraction.}
    \label{fig-entropyequiv}
\end{figure}

Based on the last property, one can think of the indicator $f_{eq}$ as an equivalent volume fraction, corresponding to a probability distribution $\hat{p}$, homogeneous in a subdomain $\Hat{\mathcal{B}}$, with the same relative entropy as $p$. Furthermore, it is possible to show that $|\Hat{\mathcal{B}}|\leq |\mathcal{B}_p|$, which means that, for equal entropies, the smallest possible support attainable corresponds to the partially homogeneous distribution. The relation between densities $p$ and $\hat{p}$ and their supports is represented in Figure~\ref{fig-entropyequiv}. 

To consider the case of ${\int_\mathcal{B} {s} \ dV = 0}$, note that the support of $s$ is also $\mathcal{B}_p$ and that integration over $\mathcal{B}_p$ instead of $\mathcal{B}$ does not affect the result in Eqs.~\eqref{eq-probdistr} and \eqref{eq-KLdivergence}. Therefore, we have the following limits
\begin{equation}
    \lim_{|\mathcal{B}_p|\to 0} D(p||q_h) = \infty \qquad \implies \qquad \lim_{|\mathcal{B}_p|\to 0} f_{eq} = 0 \ .
\end{equation}

We adopt the limit of $f_{eq} = 0$ by convention, which agrees naturally with the interpretation of $f_{eq}$ as an equivalent volume fraction: $s=0$ almost everywhere implies that there is no $s$-content at all in the volume. 

In the context of measuring the diffusion of the internal variables in an inelastic process, the results shown in this paper are based on choosing $s = \xi$ in Eq.~\eqref{eq-probdistr}, where $\xi$ was introduced in Section~\ref{sect-epmodel} to be the equivalent plastic strain. Therefore, $f_{eq}$ represents an equivalent plastified volume fraction. Note, however, that the value of this indicator may decrease as the simulation progresses if the plastic slip becomes more concentrated (for example due to the localization of the plastic deformation). Moreover, following the observation in the previous paragraph, even if the whole domain has some value of plastic slip, the indicator will not attain its maximum value unless the plastic slip is homogeneous, which reflects the possibility of further diffusion.

\section{Numerical examples}
\label{sect-examples}


This section presents examples in which deformable bodies undergo cyclic large distortions. The latter
are severe enough to justify remeshing the body during a simulation. The purpose
of the examples is to study the numerical diffusion introduced
by remeshing operations, both for the standard finite element method and for the nodal strain
approach.

In all examples, the body is modelled with the elasto-plastic constitutive law presented in Section~\ref{sect-epmodel}, and a neo-Hookean energy function such as the one in Eq.~\eqref{eq-ghat2} is adopted as the stabilizing function $\hat{g}$. The geometry, boundary conditions, and material properties are described individually for each problem. All problems are solved quasi-statically.

The meshes used in the examples consist of tetrahedral elements, as is natural in the context
of general meshing procedures. In this line, the standard finite element results correspond to a
displacement-pressure mixed formulation with a P1/P1 interpolation and a projection
stabilization \cite{bochev2004wk}. The internal variables of the standard finite element approach are
remapped following a widely used projection strategy based on element shape functions \cite{peric1996transfer}. To make the comparisons as fair as possible, remeshing operations are performed at
regular pseudo-time intervals, irrespective of the state of the meshes, and the solutions are
projected over the whole body. The time intervals for remeshing are chosen based on the total time
of the simulation to generate several remeshing events in each example.

\subsection{Torsion of a cylinder}
\label{sect-Example1}

\begin{figure}[t]
    \centering
    \begin{subfigure}{0.25\textwidth}
        \centering
        \includegraphics[width=0.7\linewidth]{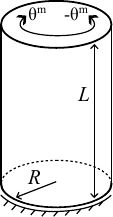}
    \end{subfigure}
    \hspace{1cm}
    \begin{subfigure}{0.25\textwidth}
        \centering
        \includegraphics[width=\linewidth]{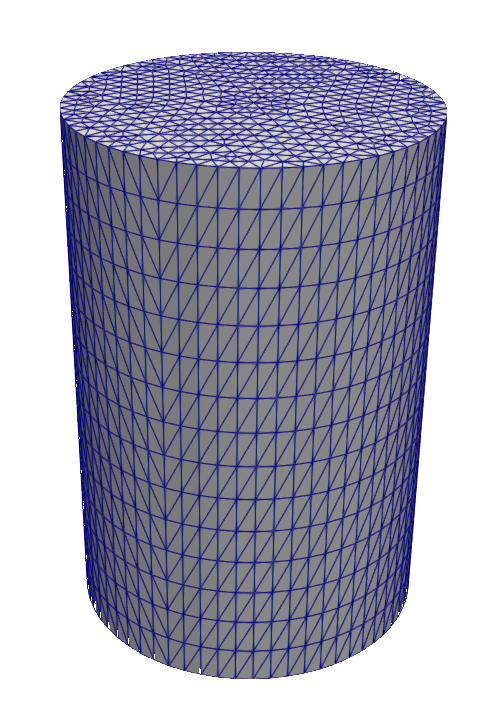}
    \end{subfigure}
    \caption{Torsion of a cylinder. Problem setting and initial mesh.}
    \label{fig-problem1}
\end{figure}

The first problem consists of a cylinder of height $L=1.5$ and radius $R=0.5$, with its bottom surface fixed
and subjected to an imposed cyclic rotation $\bar{\theta}$ on its top surface (see Figure~\ref{fig-problem1}) while fixing its vertical displacement. The imposed rotation increases linearly according to a constant angular velocity $\omega$, until it reaches a maximum rotation $\theta^m$. At that point the direction of the rotation is reversed, until the condition $\bar{\theta}=-\theta^m$ is attained. This loading history can be expressed as
\begin{equation}
    \bar{\theta}(t) := (-1)^k (\omega t -2\,k\,\theta^m) \quad \text{for} \ \ (2k-1)\theta^m/\omega < t \leq (2k+1)\theta^m/\omega \quad \text{and} \ \ k = 0,1,...
    \label{eq-load-cylinder}
\end{equation}

The direction of the rotation is reversed three times, and the simulation stops when the imposed rotation returns to zero, that is, after three half cycles. In Eq.~\eqref{eq-load-cylinder}, this corresponds to $k=0,1,2,3$ and $0\leq t\leq 6\theta^m/\omega$. The loading parameters are $\omega=1$ and $\theta^m=0.6\pi$. The (dimensionless) parameters of the material are Young's modulus $E=500$, Poisson's ratio $\nu = 0.3$, isotropic hardening modulus $K = E/20$, kinematic hardening modulus $H=0$, and yield stress $\sigma_{Y} = E/500$. The stabilizing function takes the same elastic moduli and a stabilization coefficient $\alpha = 0.2$. The body is initially meshed with 36000 elements, see Figure~\ref{fig-problem1}. Remeshing is performed every unit of pseudo-time.

\begin{figure}
    \centering
    \includegraphics[width=1\linewidth]{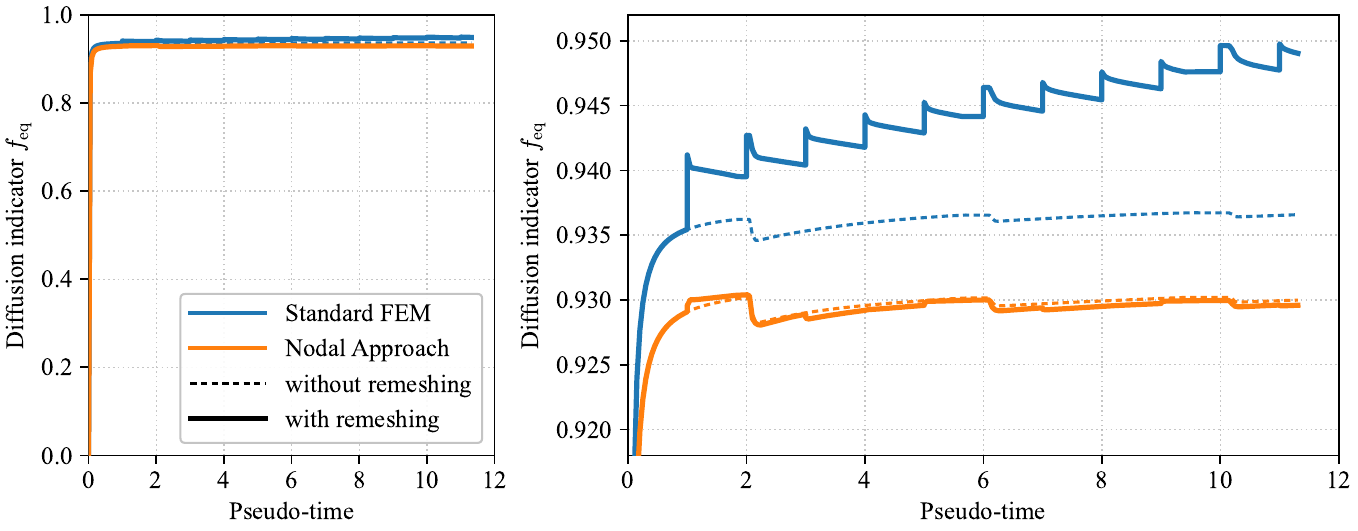}
    \caption{Internal variable diffusion due to body remeshing for the standard FEM and nodal strain approach for a cylinder subjected to an imposed rotation.}
    \label{fig-example01}
\end{figure}

Figure~\ref{fig-example01} compares the evolution of the diffusion indicator $f_{eq}$ of
Eq.~\eqref{eq-diffindicator} obtained with and without remeshing for both techniques. The general
trend --- irrespective of whether the body is remeshed or not --- is that the equivalent plastified volume fraction reaches an almost constant value at the very beginning of the deformation process. The reason is that, although the equivalent plastic deformation field $\xi$ grows with the imposed rotation, it does so in an almost proportional fashion, so that the probability distribution~\eqref{eq-probdistr} remains nearly constant in time.

More central to the scope of this paper are the discontinuities that remeshing introduces in the
curves. Both methods exhibit jumps in the equivalent plastified volume fraction each time the body
is remeshed, which reveals that a non-physical diffusion is introduced to the solution. For the standard finite element formulation, however, the magnitude of these jumps is
clearly much larger. Moreover, the numerically induced diffusion accumulates over successive
remeshing steps. 

It should be stressed that, for the purpose of this paper, the meaningful comparison in terms of the diffusion indicator comes from the \emph{jumps at remeshing events}, and not from the absolute levels attained by each method. The two curves are built from different discrete representations of the field $\xi$ --- dual-cell constants for the nodal approach and Gauss point values for the standard finite element method --- so that a systematic offset between them is expected and carries no information about the quality of either solution. The vertical separation between the solid and dashed curves of a single method, on the other hand, measures exclusively the effect of remeshing, since the two simulations are otherwise identical.

\subsection{Distortion of a block}
\label{sect-Example2}

\begin{figure}[t]
    \centering
    \begin{subfigure}{0.45\textwidth}
        \centering
        \includegraphics[width=\linewidth]{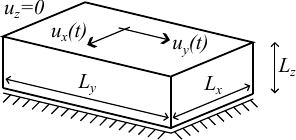}
    \end{subfigure}
    \hspace{1cm}
    \begin{subfigure}{0.45\textwidth}
        \centering
        \includegraphics[width=\linewidth]{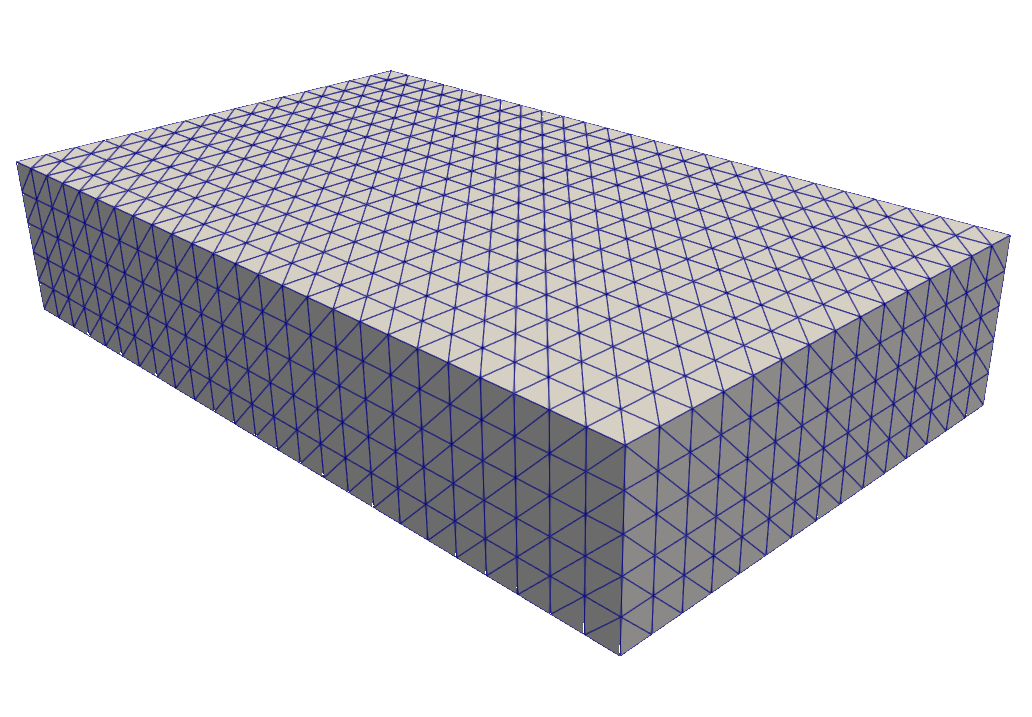}
    \end{subfigure}
    \caption{Distortion of a block. Problem setting and initial mesh.}
    \label{fig-problem2}
\end{figure}

The second problem analyses a rectangular prism of dimensions $L_x$, $L_y$, and $L_z$, with its bottom
surface fixed and time-dependent displacements imposed on its top surface along the $x$ and $y$
directions. This setting is sketched in Figure~\ref{fig-problem2}. The imposed displacements follow a sinusoidal law,
\begin{equation}
    \begin{split}
        u_x(t) := A_x \sin(2\pi \ t/T_x) \ , \\
        u_y(t) := A_y \sin(2\pi \ t/T_y) \ ,
    \end{split}
    \label{eq-dispexample2}
\end{equation}

The dimensions of the block are $L_x = 3$, $L_y=5$, and $L_z=1$, and the loading parameters are $A_x = 0.8$, $T_x = 1$,  $A_y = 1.2$, and $T_y = 1.5$. The pseudo-time interval of the complete simulation is $0 \leq t \leq 3$, so that three periods are completed in the oscillation of the $x$-displacements and two periods in that of the $y$-displacements. The real and stabilizing material parameters coincide with those of the previous example, namely Young's modulus $E=500$, Poisson's ratio $\nu = 0.3$ (both constitutive relations), isotropic hardening modulus $K = E/20$, kinematic hardening modulus $H=0$, yield stress $\sigma_{Y} = E/500$, and stabilization coefficient $\alpha = 0.2$. The body is initially meshed with 11250 elements (see Figure~\ref{fig-problem2}), and remeshing is performed every 0.25 units of pseudo-time.

\begin{figure}[t]
    \centering
    \includegraphics[width=1\linewidth]{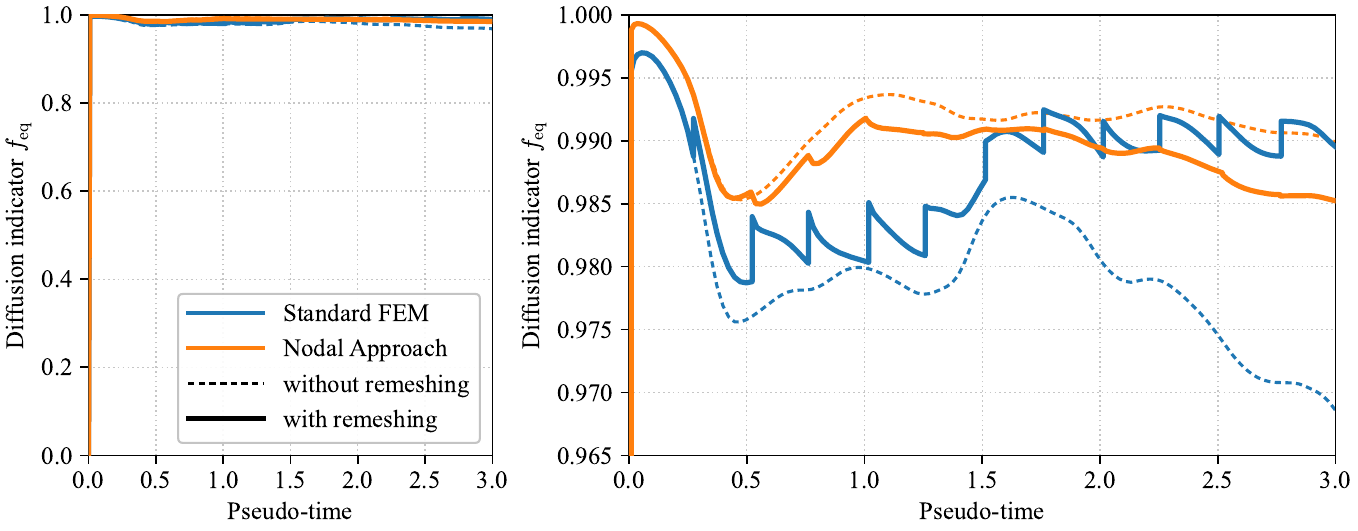}
    \caption{Internal variable diffusion due to body remeshing for the standard FEM and nodal strain approach for a block subjected to an imposed distortion.}
    \label{fig-example02}
\end{figure}

The evolution of the diffusion indicator $f_{eq}$ is compared in Figure~\ref{fig-example02}. As in the previous example, $f_{eq}$ grows rapidly at the beginning of the deformation process, and the overall behaviour is consistent with the one reported there: every remeshing event produces a noticeable jump in the diffusion indicator for the standard FEM formulation, while at the same scale of observation the curve of the nodal approach is virtually continuous.

Unlike in the example of Section~\ref{sect-Example1}, the reference curves obtained without
remeshing are not (almost completely) monotonic here: $f_{eq}$ experiences successive stages of growth
and decay. The jumps of the standard finite element curve due to the remeshing events occur at instants where the indicator is increasing, decreasing, or close to a local extremum, which indicates that the diffusion introduced by the transfer operation is governed by the remeshing operation itself rather than by the stage of the deformation process at which it takes place.

This example also reveals a feature of the solutions that involve remeshing. The evolution of $f_{eq}$
after a remeshing event differs between the curves computed with and without remeshing: an
increasing trend in the dashed line does not necessarily correspond to an increasing trend in the
solid one, and the two curves do not preserve the offset generated by the jump. This phenomenon is
observed for both methods and shows that remeshing has an effect on the solution itself. In this
sense it is worth insisting that the ability of a method to limit the numerically induced diffusion
of information is reflected in the magnitude of the jumps, and not in the absolute values of the
indicator.

\subsection{Rotation of a tube}
\label{sect-Example3}

\begin{figure}[htbp]
    \centering
    \begin{subfigure}{0.25\textwidth}
        \centering
        \includegraphics[width=\linewidth]{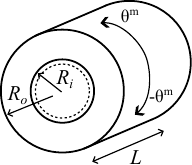}
    \end{subfigure}
    \hspace{1cm}
    \begin{subfigure}{0.22\textwidth}
        \centering
        \includegraphics[width=\linewidth]{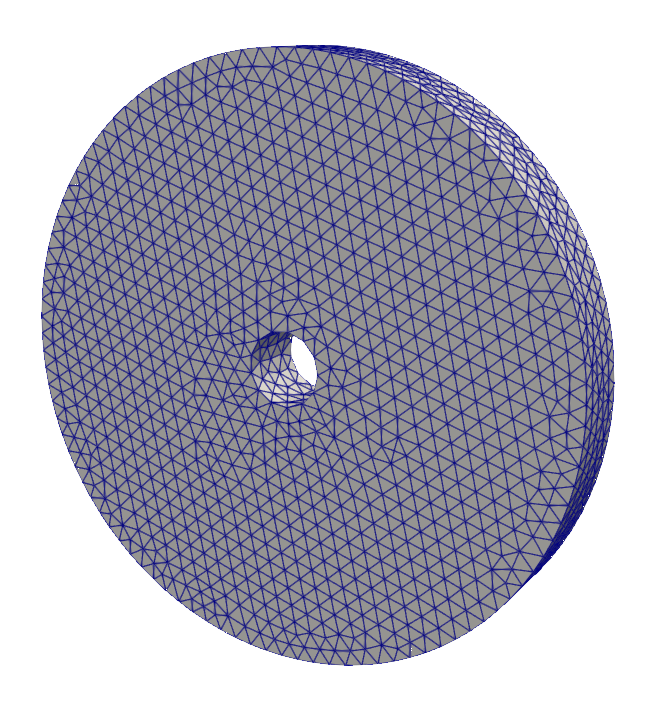}
    \end{subfigure}
    \caption{Rotation of a tube. Problem setting and initial mesh.}
    \label{fig-problem3}
\end{figure}

The last problem consists of a tube of length $L$, inner radius $R_i$, and outer radius $R_o$. Its two plane surfaces and its inner curved surface are fixed, and a cyclic rotation $\bar{\theta}$ about the longitudinal axis is imposed on the outer surface. This setting is represented in Figure~\ref{fig-problem3}. The imposed rotation follows expression~\eqref{eq-load-cylinder}, already used in Section~\ref{sect-Example1}.

In this case, the geometry and the material parameters are chosen so as to generate a more localized evolution of the plastic deformation. The purpose of this modification is to examine the influence of remeshing on the diffusion indicator when the latter takes low values. The dimensions of the tube are $L=2.0$, $R_i=1.0$, and $R_o=8.0$. The loading parameters are $\omega=1$ and $\theta^m = 0.40\pi$, the direction of the rotation is reversed five times, that is $k=0,...,5$ in Eq.~\eqref{eq-load-cylinder}, and the pseudo-time interval of the simulation is $0\leq t \leq 10\theta^m/\omega$. The parameters of the real material are Young's modulus $E=500$, Poisson's ratio $\nu = 0.3$, isotropic hardening modulus $K = E/20$, kinematic hardening modulus $H=0$, and yield stress $\sigma_{Y} = E/10$, whereas the stabilizing function takes the same elastic moduli and a stabilization coefficient $\alpha = 0.2$. The body is initially meshed with 16064 elements, see Figure~\ref{fig-problem3}. Remeshing is performed every 1.5 units of pseudo-time.

\begin{figure}
    \centering
    \includegraphics[width=0.8\linewidth]{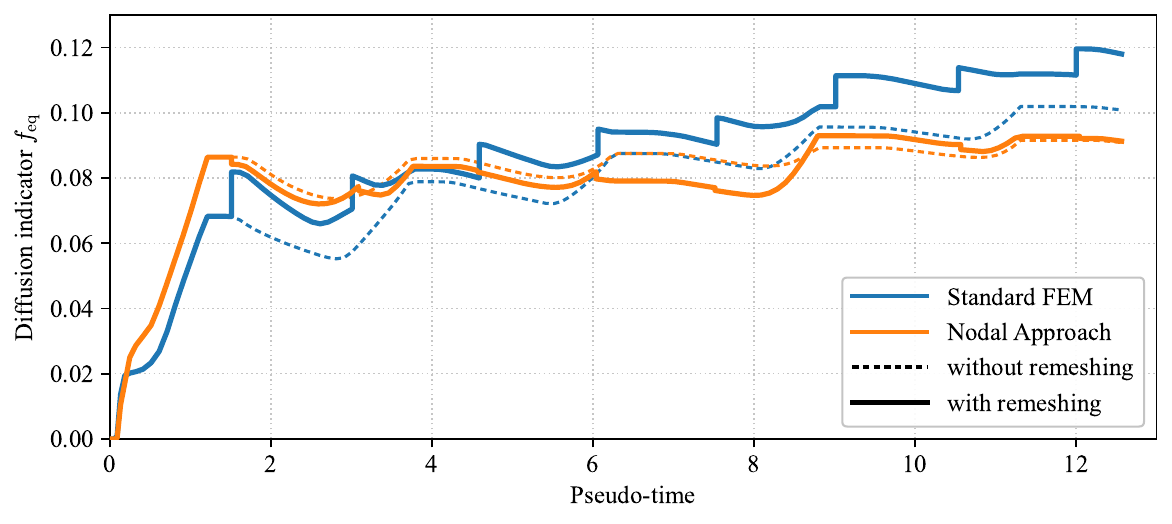}
    \caption{Internal variable diffusion due to body remeshing for the standard FEM and nodal strain approach for a tube subjected to an imposed rotation on its outer surface.}
    \label{fig-example03}
\end{figure}

The evolution of the diffusion indicator $f_{eq}$ for the four cases is depicted in Figure~\ref{fig-example03}. In contrast with the first two examples, the more localized distribution of the plastic strain keeps the indicator at low values throughout the simulation, and the reference curves themselves exhibit a marked evolution, with alternating stages of growth and stagnation associated with the successive reversals of the imposed rotation. As a consequence, the jumps of the standard FEM curve, although similar in absolute magnitude to those of the previous examples, represent a much larger fraction of the value of the indicator itself. This last observation is the main motivation for the example. Since $f_{eq}$ admits the interpretation of an equivalent plastified volume fraction, relatively large increments of the indicator correspond to a spreading of the plastic slip over a significant fraction of the total volume of the body.

\section{Conclusions}
\label{sect-conclusion}

This work has investigated a nodal strain finite element formulation for small
and large strain solid mechanics. In particular, we have followed the work of Lamichhane for small
strain elasticity~\cite{lamichhane2009hu} and later extended it to the finite strain regime, proposing the first variationally formulated, stabilized, nodal strain finite element method for this class of problems in mechanics.

The nodal strain formulation can be interpreted, for both the small- and large-strain kinematics, as a mixed
finite element method where the interpolated fields are piecewise polynomials on a primal
mesh or on a dual mesh. This split allows the elimination of the fields on the dual mesh, leading to a
single-field variational method that can be interpreted as an \emph{assumed strain} finite element,
and it is locking-free.

In the finite strain regime, the choice of stabilizing term is non-trivial and we have proposed
a general class of stabilizing functions that ensure the variational consistency of the assumed
strain formulation, a (relatively) simple implementation, and objectivity. While there is
quite some leeway in this choice, we have identified physical arguments that can be used
to select certain types of stabilization functions. Numerical examples confirm these choices.

One of the main advantages of the nodal strain formulation is that the material points are located at
the same points as the nodes of the finite element mesh. This has far-reaching implications: first,
the method can be interpreted as a true particle formulation. All the kinematic and material history
information (if any) is stored at the nodes, while Dirichlet and Neumann boundary conditions are
easily imposed. Second, if remeshing needs to be performed, there is no need to project material data
onto the material points of the new mesh, because the latter do not change --- as long as no new
nodes are introduced. If new nodes had to be added, the internal variables of the material would
need to be interpolated to the new nodes.

Avoiding the transfer of material information between different meshes reduces the diffusion of
information that causes, among other phenomena, the spreading of history variables. We have
introduced a measure that has allowed us to assess the size of this ``spreading'' and used it in
several examples. These tests show, invariably, that nodal strain methods are less afflicted by
information diffusion of internal variables and thus, together with their robustness with regard to
locking, indicate that they are well-suited for numerical simulations of finite strain mechanics
with large distortions.

\section*{Acknowledgements}
The authors acknowledge the funding received from projects 
PLEC2023-010190 and PID2025-174275NB-I00 from the Spanish Ministry of Science and Innovation.

\appendix

\section{Appendix: remeshing performance in numerical experiments}
\label{sect-appendix}

This appendix summarizes information related to the performance of the flip-based remeshing algorithm for the sake of completeness. The implementation allows for four types of interior tetrahedral flip operations (those that do not affect the boundary of the domain) $2\to3$, $3\to2$, $4\to4$, and $5\to6$ (a detailed explanation of face and edge swapping operations can be found in the literature \cite{freitag1997tetrahedral}), and one type of boundary flip denoted by $2\to2$ which involves a face swap of two tetrahedra with one face each in a certain boundary of the domain. The reported information includes the number of performed interior flips (NIF), performed number of boundary flips (NBF), and variation in the total number of elements of the mesh (VTE).

Remeshing strategies based exclusively on flip/swap operations guarantee that the quality of the mesh does not worsen. However, it is known that they do not guarantee the generation of a new mesh satisfying all the imposed \emph{good-quality tolerances}. This, however, is not a limitation for the purpose of this paper. In this sense, it is remarked that remeshing the body has two major effects: the redefinition of finite dimensional function spaces of the discrete problem and the projection of the solution (which includes the remapping of the internal variables). A greater number of flip operations, or equivalently a greater number of modified elements, has a larger impact on the change of the function spaces. However, since the projection of the solution is performed over the whole domain, irrespective of how much the mesh was actually changed, the magnitude of the diffusion of the information is rather insensitive to the configuration of the new mesh, which might be the result of a purely flip-based strategy or more sophisticated ones such as Delaunay triangulations.

Two measures were used to assess the quality of a tetrahedral element: the minimum dihedral angle, defined as the minimum angle between any pair of faces, and the mean-ratio metric, defined as a normalized ratio between the tetrahedron volume and the squared lengths of its edges. Once threshold values, $\beta_{\min}$ and $r_{\min}$ respectively, are specified, any tetrahedron whose quality metrics fall below the corresponding thresholds is classified as a poor-quality element, and the remeshing algorithm searches for improved configurations using the topological flips described above. 

An additional criterion is imposed for boundary flips, consisting of a maximum allowable change in the total volume of the affected patch, denoted by $v_{\max}$. Any candidate boundary flip that produces a relative volume change greater than this prescribed tolerance is rejected. This control is particularly important for curved boundaries, where a flip may modify the polyhedral approximation of the boundary geometry. 

The tolerances adopted in this work were $\beta_{\min}=15^\circ$, $r_{\min}=0.3$, and $v_{\max}=0.2$. Tables \ref{table-flips-example1}, \ref{table-flips-example2}, and \ref{table-flips-example3} summarize the information for the examples of Sections~\ref{sect-Example1}, \ref{sect-Example2}, and ~\ref{sect-Example3}, respectively.

\begin{table}[ht]
\centering
    \begin{tabular}{c c c c c c c}
         Example 1 & \multicolumn{3}{c}{Standard FEM (Gauss integration)} & \multicolumn{3}{c}{Nodal approach (dual mesh integration)} \\
         Remesh time & NIF & NBF & VTE & NIF & NBF & VTE \\
        \hline
         1.0 & 232 & 0 & 0 & 264 & 0 & 0 \\
         2.0 & 867 & 365 & 0 & 798 & 234 & 0 \\
         3.0 & 0 & 0 & 0 & 0 & 0 & 0 \\
         4.0 & 1 & 0 & 0 & 10 & 0 & 0 \\
         5.0 & 325 & 17 & 0 & 356 & 12 & 0 \\
         6.0 & 753 & 4 & 0 & 555 & 12 & 0 \\
         7.0 & 0 & 0 & 0 & 1 & 0 & 0 \\
         8.0 & 12 & 0 & 0 & 22 & 0 & 0 \\
         9.0 & 482 & 2 & 0 & 371 & 0 & 0 \\
         10.0 & 3 & 0 & 0 & 5 & 0 & 0 \\
         11.0 & 0 & 0 & 0 & 0 & 0 & 0 \\
    \end{tabular}
    \caption{Remeshing performance of numerical example of Section~\ref{sect-Example1}. Initial mesh with 36000 tetrahedral elements.}
    \label{table-flips-example1}
\end{table}

\begin{table}[ht]
\centering
    \begin{tabular}{c c c c c c c}
         Example 2 & \multicolumn{3}{c}{Standard FEM (Gauss integration)} & \multicolumn{3}{c}{Nodal approach (dual mesh integration)} \\
         Remesh time & NIF & NBF & VTE & NIF & NBF & VTE \\
        \hline
         0.25 & 1061 & 98 & 0 & 893 & 130 & 0 \\
         0.50 & 55 & 5 & 0 & 22 & 4 & 0 \\
         0.75 & 37 & 26 & +2 & 7 & 4 & 0 \\
         1.00 & 398 & 115 & -4 & 131 & 82 & 0 \\
         1.25 & 567 & 76 & +5 & 772 & 83 & -4 \\
         1.50 & 1442 & 58 & -3 & 249 & 121 & +3 \\
         1.75 & 1478 & 22 & +4 & 815 & 154 & -9 \\
         2.00 & 1261 & 239 & -15 & 448 & 58 & -1 \\
         2.25 & 361 & 53 & -1 & 165 & 40 & 0 \\
         2.50 & 1455 & 45 & +7 & 320 & 105 & +3 \\
         2.75 & 762 & 105 & -1 & 521 & 79 & +1 \\
    \end{tabular}
    \caption{Remeshing performance of numerical example of Section~\ref{sect-Example2}. Initial mesh with 11250 tetrahedral elements.}
    \label{table-flips-example2}
\end{table}

\begin{table}[ht]
\centering
    \begin{tabular}{c c c c c c c}
         Example 3 & \multicolumn{3}{c}{Standard FEM (Gauss integration)} & \multicolumn{3}{c}{Nodal approach (dual mesh integration)} \\
         Remesh time & NIF & NBF & VTE & NIF & NBF & VTE \\
        \hline
         1.5 & 136 & 14 & -24 & 103 & 8 & -21 \\
         3.0 & 106 & 14 & -8 & 81 & 11 & -12 \\
         4.5 & 23 & 4 & -5 & 2 & 0 & -1 \\
         6.0 & 111 & 28 & -3 & 98 & 24 & +2 \\
         7.5 & 24 & 5 & -1 & 16 & 0 & +1 \\
         9.0 & 149 & 51 & +2 & 84 & 24 & -2 \\
         10.5 & 121 & 25 & +4 & 76 & 14 & -3 \\
         12.0 & 19 & 6 & +1 & 7 & 2 & -1 
    \end{tabular}
    \caption{Remeshing performance of numerical example of Section~\ref{sect-Example3}. Initial mesh with 16064 tetrahedral elements.}
    \label{table-flips-example3}
\end{table}

\bibliographystyle{unsrt}
\bibliography{biblio}

\end{document}